\documentclass[preprint,11pt]{elsarticle}

\usepackage{graphicx}
\usepackage{verbatim}
\usepackage{amssymb}
\usepackage{lineno}
\usepackage{amsthm}
\usepackage{graphicx}
\usepackage{caption}
\usepackage{subcaption}
\usepackage[margin=1in]{geometry}
\usepackage{booktabs}
\usepackage{tabularx}
\usepackage{amsmath}
\usepackage{dsfont}
\usepackage[multiple]{footmisc}
\usepackage{mathtools}
\usepackage{float}
\usepackage{multirow}
\usepackage{natbib}
\usepackage{color}
\biboptions{sort&compress}
\restylefloat{table}
\usepackage{upgreek}
\usepackage{nicefrac}
\usepackage{xargs}                      
\usepackage[pdftex,dvipsnames]{xcolor} 
\usepackage[colorinlistoftodos,prependcaption,textsize=tiny]{todonotes}
\usepackage{xcolor}
\usepackage{multicol}
\usepackage{nomencl}
\makenomenclature
\renewcommand{\vec}[1]{\mathbf{#1}}

\theoremstyle{plain}

\theoremstyle{definition}

\theoremstyle{remark}

\newcommandx{\unsure}[2][1=]{\todo[linecolor=red,backgroundcolor=red!25,bordercolor=red,#1]{#2}}
\newcommandx{\change}[2][1=]{\todo[linecolor=red,backgroundcolor=red!25,bordercolor=red,#1]{#2}}
\newcommandx{\info}[2][1=]{\todo[linecolor=OliveGreen,backgroundcolor=OliveGreen!25,bordercolor=OliveGreen,#1]{#2}}
\newcommandx{\improvement}[2][1=]{\todo[linecolor=Plum,backgroundcolor=Plum!25,bordercolor=Plum,#1]{#2}}
\newcommandx{\thiswillnotshow}[2][1=]{\todo[disable,#1]{#2}}

\usepackage{lineno}

\journal{Computer Methods in Applied Mechanics and Engineering}

\begin{document}

\begin{frontmatter}


\title{Data-driven surrogates for high dimensional models using Gaussian process regression on the Grassmann manifold}



\author{D. G. Giovanis$^{a}$, M. D. Shields$^{a}$}

\address{$^{a}$Department of Civil \& Systems Engineering,
	Johns Hopkins University, Baltimore, MD, 21218, USA \\}

\begin{abstract}

\noindent
This paper introduces a surrogate modeling scheme based on Grassmannian manifold learning  to be used for cost-efficient predictions of  high-dimensional stochastic systems. The method exploits subspace-structured features of each solution by projecting it onto a Grassmann manifold. This point-wise linear dimensionality reduction harnesses the structural information to assess the similarity between solutions at different points in the input parameter space. The method utilizes a solution clustering approach in order to identify regions of the parameter space over which solutions are sufficiently similarly such that they can be interpolated on the Grassmannian.  In this clustering, the reduced-order solutions are partitioned into disjoint clusters on the Grassmann manifold using the eigen-structure of properly defined Grassmannian kernels and, the Karcher mean of each cluster is estimated.  Then, the points in each cluster are projected onto the tangent space with origin at the corresponding Karcher mean using the exponential mapping. For each cluster, a Gaussian process regression model is trained that maps the input parameters of the system to the reduced solution points of the corresponding cluster projected onto the tangent space. Using this Gaussian process model, the full-field solution can be efficiently predicted at any new point in the parameter space. In certain cases, the solution clusters will span disjoint regions of the parameter space. In such cases, for each of the solution clusters we utilize a second, density-based spatial clustering to group their corresponding input parameter points in the Euclidean space. The proposed method is applied to two numberical examples. The first is a nonlinear stochastic ordinary differential equation with uncertain initial conditions where the surrogate is used to predict the time history solution. The second involves modeling of plastic deformation in a model amorphous solid using the Shear Transformation Zone theory of plasticity, where the proposed surrogate is used to predict the full strain field of a material specimen under large shear strains.

\end{abstract}

\begin{keyword}
Grassmann manifold \sep Spectral clustering \sep Gaussian process regression \sep Machine learning \sep Nonlinear projection \sep Interpolation


\end{keyword}

\end{frontmatter}


\section{Introduction}
\label{S:Intro}

The need for very high-fidelity, computationally expensive models to accurately capture the underlying physics of engineering systems is a limiting factor in computational mechanics. Despite the huge strides made over the last decades to solve these massive-scale problems in engineering (i.e. higher-order mathematical models, better approximation methods, advances towards exascale computing and  development of sophisticated numerical  algorithms), large-scale modeling of real-world applications remains an uphill battle, especially in the framework of uncertainty quantification (UQ) where a large number of repeated simulations are required. However, efforts to overcome these computational barriers are continuous. 

A drastic reduction in the computational cost  can be achieved  by developing a simpler mathematical function (surrogate) of the model that is inexpensive to evaluate and is able to reproduce the system's input-output relation. In order to build a surrogate that can approximate/interpolate the solution between sample points where the exact solution is known, a finite number of high-fidelity model evaluations are required at known points in the input parameter space. Among the most widely-used surrogate modeling techniques in the framework of UQ,  include generalized polynomial chaos expansions (PCE) \cite{xiu2002wiener,MEgPC2005, MEgCM2006, Blatman2011, papaioannou2019pls, luthen2020sparse}, which are also widely-used for the related stochastic collocation methods \cite{Babuska2007,MEPCM2008, ASGC2009, Xiu2005, loukrezis2019assessing}) and the Gaussian process (GP) regression model \cite{Krige51, Santner2003, Rasmussen2006, Echard2011, schobi2015polynomial} which is a machine learning technique constructed by fitting a Gaussian random function to the training data. However, the training time  of a GP  scales cubically  $O(n^3)$ with the number of available training data $n$ while the memory scaling is $O(n^2)$. To this end, efforts to overcome scaling problems are made with the propose of sparse GP approximations \cite{Csato2002, Smola2001, wang2020accelerated}.

An alternative method to reduce computational cost is through reduced-order modelling. Reduced-order models (ROMs) seek the simplest mathematical representation of the high-fidelity model that captures the  dominant behavior of the system, i.e. aim to preserve the governing equations while drastically reducing the number of degrees of freedom.  ROM approaches include reduced-basis methods \cite{Rozza2008, Tan2008} and nonlinear projection methods \cite{Benner2015}. From the first category,  a widely-used approach employs the proper orthogonal decomposition (POD) \cite{POD2000}. The key idea of POD is that the subspace spanned by a number of system solutions constitutes a basis of the total solution space that can accurately describe the system. Recent efforts in UQ have focused on adaptive methods to identify the reduced basis, or local reduced bases (see e.g. \cite{zou2019adaptive}) as well as Bayesian approaches to improve the identification of the reduced basis \cite{stabile2020bayesian}. On the other hand, nonlinear projection methods are based on trajectories between solutions computed using high-dimensional models that are contained and may be interpolated in low-dimensional subspaces.  Related work can be found in \cite{Amsallem2008, Amsallem2011} in which a groundbreaking  nonlinear projection method for interpolating ROMs on the Grassmann manifold was introduced. 

Additionally, recent advances in machine learning have forged a new path in utilizing manifold learning for processing high-dimensional data and thus, enhancing the quality of the surrogate models in the framework of data-driven UQ.  Manifold learning is built on the assumption that  the high-dimensional data lie on a lower-dimensional space (manifold) which encompasses the structural information of the data. Identifying this manifold can lead to lower complexity, better insight into the physics of the problem and better computational performance.  Over the last decade we have witnessed the emergence of various linear dimension reduction methods such as POD/principal component analysis (PCA) \cite{PCA1986} as well as non-linear dimension reduction methods, such as locally linear embedding (LLE) \cite{Roweis2000}, Laplacian eigenmaps \cite{Belkin2003}, Diffusion maps (DMAPs) \cite{DMaps2005},  ISOMAP \cite{Tenenbaum2000}, Hessian eigenmaps \cite{Dohoho2003}, Kernel PCA \cite{Scholkopf1998} and local tangent space alignment \cite{Zhang2004}.  These methods for dimension reduction are now being fully integrated into surrogate-based uncertainty quantification efforts \cite{lataniotis2020extending}. 

In this work, we propose a method that combines manifold learning with surrogate model construction for interpolation of reduced-order solutions that can be used, through inversion of the dimensionality reduction, to predict the full solution at new points in the parameter space without requiring expensive model evaluation. Strictly speaking, the model we construct is a surrogate model and not a ROM because we do not solve the governing physical/mathematical equations in the reduced basis. However, it can be viewed as a ``physics-informed'' surrogate model in the sense that we interpolate directly in the reduced space of physically/mathematically relevant bases for the specific problem -- although it is not physics-informed in the sense that the solutions are constrained to satisfy physical principals (e.g. conservation of mass, energy, momentum, etc.) \cite{raissi2019physics}. 

More specifically, we first evaluate the expensive high-dimensional computational model at a set of training points. At each training point, we reduce the dimension of the full solution through a projection onto the Grassmann manifold. We then exploit the structure of the points on the manifold to identify clusters of points having similar solutions. For each cluster, a local GP is constructed (in the tangent space of the manifold) that creates a mapping between the input parameter space and a low-dimensional (manifold) representation of the high-dimensional solution. Straightforward inversion of the projection then allows us to approximately reconstruct the high-dimensional solution at any arbitrary point in the input parameter space without evaluating the expensive computational model. Applications are considered that involve time-dependent having potentially long duration and spatially varying solutions with a very large number of degrees of freedom.


It is worth mentioning that other works to this, and related ends have been conducted very recently.  Wang et al.\ \cite{wang2020probabilistic} use manifold learning methods to identify the intrinsic structure of particular solutions (referred to as probabilistic performance patterns) and how they relate to inputs to the system. Kalogeris and Papadopoulos have presented a new surrogate modeling approach for high dimensional models that leverages a nonlinear dimension reduction through diffusion maps \cite{kalogeris2019diffusion}. Soize and Ghanem \cite{soize2016data} developed a manifold learning technique to identify the underlying probability structure of high-dimensional stochastic data on a diffusion manifold and propose a novel method to sample from this distribution on the manifold. Additionally, in a series of recent papers Soize, Farhat, and collaborators have established a manifold projection-based framework for learning model-form uncertainties \cite{soize2017nonparametric,farhat2018modeling,soize2019probabilistic}.  In two prior works \cite{Giovanis2018, Giovanis2019}, the authors utilize an unstructured multi-element mesh of the parameter space in combination with polynomial chaos expansion-based surrogates to interpolate high dimensional solutions using a locally reduced basis \cite{Giovanis2019} and the Grassmann manifold projection of the solution \cite{Giovanis2018}. The primary shortcoming of these methods relate to the requisite multi-element discretization of the parameter space, which limits the approach to problems with low-dimensional input parameters. The proposed approach does not suffer from the same restrictions.

Finally, an important emerging area of application of combined ML and UQ methods is in materials modeling at various scales ranging from atomistic simulations \cite{wang2019modeling} to continuum constitutive modeling \cite{chernatynskiy2013uncertainty, wang2020uncertainty}. Here, we present the proposed method in the setting of continuum modeling of localized plastic deformation in amorphous solids using the Shear Transforamtion Zone theory of plasticity \cite{Falk1998, Bouchbinder2009c} as informed by molecular dynamics simulations.



The paper is organized as follows. Section 2 provides an overview of the Grassmann manifold and related concepts.  This Section introduces notions from differential geometry and thus, an elementary background in Riemannian geometry is a prerequisite.  The reader is referred to \cite{Absil2004, Edelman1998} for further reading. In Section 3, a brief review of Gaussian process regression is provided before demonstrating how it can be used to interpolate high-dimensional data on the Grassmann manifold in Section 4.  Section 5 highlights the accuracy, efficiency, and robustness of the proposed surrogate for two problems: 1) A highly nonlinear three-dimensional system of stochastic ordinary differential equations, namely the Kraichnan-Orszag (KO) three mode problem, and 2) A material modeling problem for a amorphous solids utilizing the shear transformation zone (STZ) theory of plasticity.

\section{The Grassmann manifold: An overview}
\label{S:Gr}

\noindent
Two manifolds (topological spaces that are locally similar to Euclidean space and have a globally defined differential structure) of special interest  that arise in numerical linear algebra and are considered quotients\footnote{A quotient of a group results from equivalence relations between points in the group. For example, let $\mathcal{M}_2(2\times 4)$ be the set of all $2 \times 4$ matrices of rank 2.  Two points $A, B \in \mathcal{M}_2(2\times 4)$ are called equivalent if there exists a $2 \times 2$ matrix $C$ such as $A = C\cdot B$ and $\text{det}(C)>0$.} of the special orthogonal group\footnote{$\mathcal{SO}_p$ is the smooth differential manifold (Lie group structure) of all $p \times p$ orthogonal matrices with determinant +1 (subspace of the orthogonal group ($\mathcal{O}_n$) ).} $\mathcal{SO}_p$ are the  compact\footnote{The non-compact Stiefel manifold is the set of $n \times p $ matrices that have full rank.} Stiefel manifold $\mathcal{V}_{p, n}$ and the Grassmann manifold $\mathcal{G}_{p, n}$ \cite{Absil2004,Edelman1998}.  The Stiefel manifold $\mathcal{V}_{p, n}$ is defined as the space of all $p$-dimensional orthonormal (Euclidean norm equal to 1) bases of $\mathbb{R}^{n}$, i.e the space of orthonormal matrices $\textbf{X}\in \mathbb{R}^{n\times p}$ with real entries
\begin{equation}\label{CompactStiefel}
\mathcal{V}_{p, n} = \{\textbf{X}\in \mathbb{R}^{n\times p}:\textbf{X}^\intercal\textbf{X} = \textbf{I}_p\}
\end{equation}
where $\textbf{I}_p$ is the $p\times p$ identity matrix.  

On the Grassmann manifold  $\mathcal{G}_{p, n}$ a point $\mathcal{X}$ is a linear subspace, specified by an orthogonal basis ($n\times p$ matrix). However, unlike the Stiefel manifold,  the choice of basis for the subspace is not unique. This allows interpretation of each point   on the Grassmann manifold as an equivalence of points on the Stiefel manifold (all orthonormal matrices that span the same subspace are considered equivalent). Thus, a point $\mathcal{X} \in \mathcal{G}_{p, n}$, can be represented by an orthonormal matrix $\vec{X}\in \mathbb{R}^{n \times p}$ whose columns span the corresponding subspace (Stiefel representation of the Grassmann manifold).
\begin{equation}\label{Grassmannl}
\mathcal{G}_{p, n} = \{\text{span}(\textbf{X}): \mathbf{X}\in \mathbb{R}^{n\times p}:\textbf{X}^\intercal\textbf{X} = \textbf{I}_p\}.
\end{equation}

The Grassmann manifold has a Reimannian structure, which means that it is both smooth (operations of calculus can be performed on it) and is equipped with an inner product on its tangent space. The tangent space at a point is the plane tangent to the manifold at that point, with origin at the point of tangency (the normal space is the orthogonal complement). More formally, the tangent space $\mathcal{T}_{\mathcal{X}}$ with origin point $\mathcal{X}$ is  defined by all matrices $\mathbf{\Gamma}\in \mathbb{R}^{n\times p}$ such that $\vec{X}^\intercal\vec{\Gamma}=\vec{0} $. The tangent space is a flat inner-product space  and thus, for  any two points $ \vec{\Gamma}_j, \vec{\Gamma}_j \in \mathcal{T}_{\mathcal{X}}$ we can define the inner product $\langle  \vec{\Gamma}_i, \vec{\Gamma}_j \rangle$ as the trace of the product $\vec{\Gamma}_i^\intercal \vec{\Gamma}_j$, i.e. $\langle  \vec{\Gamma}_i, \vec{\Gamma}_j \rangle = \mbox{tr}(\vec{\Gamma}_i^\intercal \vec{\Gamma}_j)$. Existence of this canonical metric makes $\mathcal{G}_{p, n}$ a Riemannian manifold and thus, lengths of paths can be defined such as the  Grassmannian geodesic which we describe next.

\subsection{Grassmannian geodesic}

\begin{figure}[!htb]
	\centering
	\includegraphics[width=0.6\textwidth]{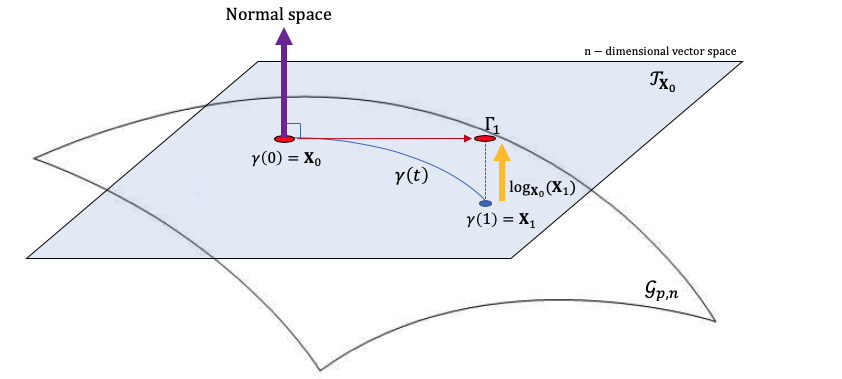}
	\caption{Illustration of tangent spaces, geodesics, and the logarithmic mapping on a manifold. Points $\mathcal{X}_0$ and $\mathcal{X}_1$ lie on the manifold with tangent space $\mathcal{T}_{\mathcal{X}_0}$. Geodesics paths are minimum distance curves on the manifold which project to straight lines in the tangent space. The logarithmic map is used to project points from the manifold onto the tangent space.}
	\label{Fig.1}
\end{figure}

\noindent
Given two orthonormal matrices $\vec{X}_0$ and $\vec{X}_1$ representing points $\mathcal{X}_0$ and $\mathcal{X}_1$ on $\mathcal{G}_{p, n}$,  respectively, the geodesic refers to the shortest curve on the manifold connecting the two points (see Fig. \ref{Fig.1}). For our purposes, geodesics are important because they provide convenient paths along which to interpolate high-dimensional objects. The geodesic  can be parameterized by a function\footnote{The function is differentiable everywhere in [0, 1]}  $\gamma(z): [0, 1] \rightarrow \mathcal{G}_{p, n}$, where  $\gamma(0)= \vec{X}_0$ and $\gamma(1)= \vec{X}_1$. The geodesic is determined by solving a second order differential equation on $\mathcal{G}_{p, n}$ that corresponds to motion along the curve having constant tangential velocity (i.e. all acceleration is normal to the manifold). Hence, it can be uniquely defined by its initial conditions $\gamma(0)=\vec{X}_0$ and $\dot{\gamma}(0)=\dot{\vec{X}}_0=\vec{\Gamma}_0\in\mathcal{T}_{\mathcal{X}_0}$. 

Practically speaking, to identify the geodesic requires us to operate in the tangent space. In particular, it can be shown that the geodesic between two points $\vec{X}_0$ and $\vec{X}_1$ on $\mathcal{G}_{p, n}$ maps to a straight line between their projections $\boldsymbol{\Gamma}_0$ and $\boldsymbol{\Gamma}_1$ in the tangent space (a geodesic is the generalization of a straight line on a Reimannian manifold) as illustrated in Figure \ref{Fig.1}. Because that the tangent space in a flat, inner-product space, it is relatively straightforward to define this straight line. First, we need to define a mapping from a point on the manifold $\vec{X}_1$ onto a point $\vec{\Gamma}_1$ in the tangent space $\mathcal{T}_{\mathcal{X}_0}$. This mapping is referred to as the logarithmic mapping. We begin by defining its inverse, the exponential mapping denoted by $\log_{\mathcal{X}_0}(\boldsymbol{\Gamma}_1)=\vec{X}_1$. Taking the thin Singular Value Decomposition (SVD) of $\vec{\Gamma}_1$, $\boldsymbol{\Gamma}_1 = \vec{U}\boldsymbol{\Sigma}\vec{V}^\intercal$, allows us to write the exponential mapping as
\begin{equation}\label{eq:4}
\vec{X}_1 = \exp_{\mathcal{X}_0}(\vec{U}\boldsymbol{\Sigma}\vec{V}^\intercal)=\vec{X}_0\vec{V}\cos(\boldsymbol{\Sigma})\textbf{Q}^\intercal + \textbf{U}\sin(\boldsymbol{\Sigma})\textbf{Q}^\intercal
\end{equation}
where $\vec{U}$ and $\vec{V}$ are orthonormal matrices, $\boldsymbol{\Sigma}$ is diagonal matrix with positive real entries, and $\textbf{Q}$ is an orthogonal $n\times n$ matrix. Requiring that $\boldsymbol{\Gamma}$ lies on $\mathcal{T}_{\mathcal{X}_0}$, defines the following set of equations
\begin{subequations}
    \begin{equation}\label{eq:6a}
    \textbf{V}\cos(\boldsymbol{\Sigma})\textbf{Q}^\intercal = \vec{X}_0^\intercal \vec{X}_1
    \end{equation}
    \begin{equation}\label{eq:6b}
    \textbf{U}\sin(\boldsymbol{\Sigma})\textbf{Q}^\intercal =\vec{X}_1- \vec{X}_0\vec{X}_0^\intercal\vec{X}_1. 
    \end{equation}
\end{subequations}
Multiplying \eqref{eq:6a} by the inverse of \eqref{eq:6b} yields
\begin{equation}\label{eq:7}
    \textbf{U}\tan(\boldsymbol{\Sigma})\textbf{V}^\intercal = (\vec{X}_1- \vec{X}_0\vec{X}_0^\intercal\vec{X}_1) (\vec{X}_0^\intercal \vec{X}_1)^{-1}
\end{equation}
Thus, the exponential mapping can be performed by taking the thin SVD of the matrix  $\textbf{M} = (\vec{X}_1- \vec{X}_0\vec{X}_0^\intercal\vec{X}_1) (\vec{X}_0^\intercal \vec{X}_1)^{-1}=\textbf{U}\boldsymbol{\Sigma}\textbf{V}^\intercal$. Finally, the logarithmic map can be defined by inversion of Eq.\ \eqref{eq:7} as \cite{Begelfor2006}:
\begin{equation}\label{eq:8}
    \log_{\mathcal{X}_0}(\vec{X}_1)=\boldsymbol{\Gamma}_1 =\textbf{U}\tan^{-1}(\boldsymbol{\Sigma})\textbf{V}^\intercal
\end{equation}

Interpolation between $\boldsymbol{\Gamma}_0$ and $\boldsymbol{\Gamma}_1$ (i.e. defining the straight line connecting $\boldsymbol{\Gamma}_0$ and $\boldsymbol{\Gamma}_1$) can be performed in the usual way, as we will discuss more later. Projecting this interpolation onto the manifold and parameterizing on $z\in[0,1]$ allows us to express the geodesic as \cite{Absil2004,Edelman1998}
\begin{equation}\label{eq:9}
    \gamma(z) = \mbox{span}\left[\left(\vec{X}_0\textbf{V} \cos (z\boldsymbol{\Sigma}) + \textbf{U} \sin ( z\boldsymbol{\Sigma})\right)\textbf{V}^\intercal\right]
\end{equation}
In Eq.(\ref{eq:9}), the rightmost $\textbf{V}^\intercal$ can be omitted and still represent the same equivalence class as $\gamma(z)$.  However, due to consistency conditions the tangent vectors used for computations must be altered in the same way and for this reason everything is multiplied by $\textbf{V}^\intercal$  \cite{Edelman1998}.

\subsection{Grassmannian distance}

\noindent
The tangent space is only defined locally. Consequently, projection of points other than the point of tangency onto the tangent space introduces some approximations. To alleviate errors associated with these projections, it is important that we map only points that are ``close'' to the point of tangency. To define ``close'' we require a measure of distance \footnote{Such a measure must be invariant under rotation $d_{\mathcal{G}_{p, n}}(\vec{X}_1, \vec{X}_2) = d_{\mathcal{G}_{p, n}}(\vec{X}_1\vec{R}_1, \vec{X}_2\vec{R}_2), \forall \vec{R}_1, \vec{R}_2 \in \mathcal{O}(n)$.} on the manifold.  However, the geodesic path between the two points $\mathcal{X}_1$ and $\mathcal{X}_2$ on $\mathcal{G}_{p, n}$ does not provide any information about how far apart the subspaces are. 

Several measures of distance on the Grassmann manifold have been proposed in the literature \cite{Hamm:2008,YeLim2014}. In general, any such distance (expressed $d_{\mathcal{G}_{p, n}}(\vec{X}_1, \vec{X}_2)$) is computed as a non-linear function of the principal angles $\theta_i$ between the two subspaces using the full SVD of the product matrix $\vec{X}_0^\intercal \vec{X}_1 = \textbf{U}\boldsymbol{\Sigma} \textbf{V}^\intercal$. The principal angles\footnote{The principal angles are bounded between $0\leq \theta_1 \leq \ldots \theta_p \leq \pi/2$} are given by $\theta_i = \cos^{-1}\sigma_i$ for $i=1,\ldots,p$ where $\sigma_i$ denote the singular values from $\boldsymbol{\Sigma} = \mbox{diag}(\sigma_1, \ldots, \sigma_p) \in \mathbb{R}^{p \times p}$. Definitions of some commonly used distances between subspaces are given in Table \ref{tab:1} \cite{Hamm:2008}. 

An issue that often arises in computing Grassmannian distances between subspaces is that the subspace may have different dimensions. Let $p, k, n \in \mathbb{N}$ be such that  $p \leq k\leq n$. For any $\vec{X}_1 \in \mathcal{G}_{p, n}$ and $\vec{X}_2 \in \mathcal{G}_{k, n}$ we need to define a distance $\delta(\vec{X}_1, \vec{X}_2)$  between a subspace $\vec{X}_1$ of dimension $p$ and a subspace  $\vec{X}_2$ of dimension $k$. This distance can equivalently be expresssed as the distance between the $p$-plane $\vec{X}_1$ and the closest $p$-plane $\vec{A}$ contained in $\vec{X}_2$ measured within $\mathcal{G}_{p, n}$ or as the distance between the $k$-plane $\vec{X}_2$ and the closest $k$-plane $\vec{B}$ containing $\vec{X}_1$ measured within $\mathcal{G}_{k, n}$ \cite{YeLim2014}.  This is expressed as $\delta(\vec{X}_1, \vec{X}_2) = d_{\mathcal{G}_{p,n}}(\vec{X}_1,\Omega_{-}(\vec{X}_2)) = d_{\mathcal{G}_{k,n}}(\vec{X}_2,\Omega_{+}(\vec{X}_1))$ where  $\Omega_{+}(\vec{X}_1):=\{\vec{B}\in \mathcal{G}_{k, n}:\vec{X}_1 \subseteq \vec{B}\}$ and $ \Omega_{-}(\vec{X}_2):=\{\vec{A}\in \mathcal{G}_{p, n}:\vec{A} \subseteq \vec{X}_2\}$ are the so-called Schubert varieties. We observe that $\delta(\vec{X}_1, \vec{X}_2)$  is a distance in the sense of a distance from a point to a set, but not a distance in the sense of a metric on the set of all subspaces of all dimensions. Defining the infinite Grassmannian $\mathcal{G}_{p, \infty}$ of $p$-planes as $\mathcal{G}_{p, \infty}= \bigcup_{n = p}^{\infty} \mathcal{G}_{p, n}$, where $\mathcal{G}_{p, n}\subset \mathcal{G}_{p, n+1} $  and the doubly infinite Grassmannian  as the disjoint union of all $p$-dimensional subspaces over all $p \in \mathbb{N}$ , $\mathcal{G}_{\infty, \infty}= \coprod_{ p=1}^{\infty} \mathcal{G}_{p, \infty}$,  Ye and Lim \cite{YeLim2014}  proposed an elegant set of distance metrics on $\mathcal{G}_{\infty, \infty}$ that capture the geometry of $\mathcal{G}_{p, n}\hspace{3pt}\forall p<n$ and generalize to the common distances between subspaces of equal dimension. However, since the doubly infinite Grassmannian is non-Hausdorff and therefore non-metrizable, one way to define a metric between any pair of subspaces of arbitrary dimensions is to  transform $\delta(\vec{X}_1, \vec{X}_2) $ into a metric on $\mathcal{G}(\infty, \infty)$. To this purpose,  Ye and Lim \cite{YeLim2014} proposed to to complete $\vec{X}_1$ to an $k$-dimensional subspace of $\mathbb{R}^n$, by adding $k-p$ vectors orthonormal to the subspace $\vec{X}_2$. Distance-wise, this is equivalent to setting  $\theta_{p+1} = \ldots = \theta_{k} = \pi/2$ (see Table \ref{tab:1}).
\begin{table}[!ht]
	\centering
	\captionof{table}{Distances on $\mathcal{G}_{p, n}$ and  $\mathcal{G}(\infty,\infty)$ in terms of principal angles.} 
	\begin{tabular}{c c c}   \toprule
		\textbf{Name} & $d_{\mathcal{G}_{p, n}}(\vec{X}_1, \vec{X}_2)$ & $d_{\mathcal{G}_{\infty,\infty}}(\vec{X}_1\vec{X}_2)$  \\\midrule
		Grassmann         & $\bigg( \sum_{i=1}^{p} \theta_i^{2} \bigg)^{\frac{1}{2}}$                                           & $\bigg(|k-p|\pi^2/4 + \sum_{i=1}^{\min(p, k)} \theta_i^2 \bigg)^{1/2}$      \\ 
		Procrustes     & $2\bigg(\sum_{i=1}^p\sin^2(\theta_i/2)\bigg)^{1/2}$         & $\bigg(|k-p| + \sum_{i=1}^{\min(p, k)} \sin^2(\theta_i/2) \bigg)^{1/2}$    \\ 
		Projection     & $\bigg(\sum_{i=1}^{p} \sin^{2}\theta_i\bigg)^{1/2}$        & $\bigg(|k-p| + \sum_{i=1}^{\min(p, k)} \sin^2\theta_i \bigg)^{1/2}$    \\  \bottomrule
		\hline
		\label{tab:1}
	\end{tabular}
\end{table}

\subsection{Karcher mean}

\noindent
Finally, given that our objective is to define an interpolant between a set of points on the manifold, we further aim to identify the point on the Grassmannian whose distance to the set of points in minimized. Identifying this point allows us to map the set of points onto the tangent space with minimum error. 

We can identify this point by identifying the Riemannian center of mass of the points on the manifold, or the so-called Karcher mean{\color{blue}\footnote{Alternatives subspace averages that can be found in literature are the extrinsic manifold mean, the $L_2$-median and the flag mean \cite{Marrinan2014}.}}, defined as the point $\mathcal{Y}$ that minimizes locally the cost function $\lambda:\mathcal{G}_{p, n} \rightarrow \mathbb{R}_{\geq 0}$ given by \cite{Karcher1977}:
\begin{equation}
    \lambda(\mathcal{Y}) = \int_{\mathcal{G}_{p, n}} d_{\mathcal{G}_{p, n}}^{2}(\mathcal{Y}, \mathcal{X})\text{d}P(\mathcal{X}))
\end{equation}
where  $\text{d}P(\mathcal{X}) = \rho(\mathcal{X})\text{d}\mathcal{G}_{p, n}(\mathcal{X})$ is a probability measure over an infinitesimal volume element $d\mathcal{G}_{p, n}(\mathcal{X})$ with  probability density function $\rho(\mathcal{X})$.  The value of the function $\lambda$ at the Karcher mean is called the Karcher variance and provides a measure of the spread of the points around the Karcher mean. For of a set of independent sample points $\{\mathcal{X}_i\}_{i=1}^{N} \subset \mathcal{G}_{p, n}$ the sample Karcher mean  is defined as the local minimizer of the
\begin{equation}
\lambda(\mathcal{Y}) =  \frac{1}{N}\sum_{i=1}^{N} d_{\mathcal{G}_{p, n}}^{2}(\mathcal{X}_i, \mathcal{Y})
\end{equation}
Various iterative algorithm such as Newton’s method or first-order gradient
descent \cite{Absil2004, Begelfor2006} can be used in order to find the sample Karcher mean, which will be denoted by a bold $\textbf{m}$. The iterative algorithm found in \cite{Turaga2012CorrectionSC} is used here.  However, a unique optimal solution is not always guaranteed \cite{Karcher1977,Begelfor2006}. As shown by Begelfor and Werman in \cite{Begelfor2006}, if the ensemble of Grassmann points has a radius of greater than $\pi/4$, then the exponential and logarithmic maps are no longer bijective, and the Karcher mean is no longer unique.

\section{Gaussian process regression}
\label{S:GP}

\noindent
Gaussian process (GP) regression (or Kriging when strictly used for interpolation) is a widely used, nonparametric, Bayesian approach used in supervised learning problems i.e. learning of input-output mappings from a small training dataset. One of the advantages of GP regression is that the predictor provides a natural measure of uncertainty in the prediction. From a mathematical point of view, a GP  is a collection of random variables with a joint Gaussian distribution that  can be completely specified by a mean function and a positive definite covariance function. 

For our purposes, consider a numerical model ($\mathcal{M}$) that is developed to solve the governing equations of some physical system, e.g.
\begin{subequations}\label{eq:setPDE}
	\begin{equation}
	\mathcal{L}(\mathbf{x},t,\boldsymbol{\xi};w((\mathbf{x},t,\boldsymbol{\xi})) = f(\mathbf{x},t,\boldsymbol{\xi}), \quad \forall \mathbf{x}\in \mathcal{D},t \in T, \boldsymbol{\xi} \in \Omega
	\end{equation}
	\begin{equation}
	\mathcal{B}(\mathbf{x},t,\boldsymbol{\xi};w((\mathbf{x},t,\boldsymbol{\xi})) = g(\mathbf{x},t, \boldsymbol{\xi}), \quad \forall \mathbf{x}\in \partial\mathcal{D}, t \in T, \boldsymbol{\xi} \in \Omega 
	\end{equation}
\end{subequations}
where $\mathcal{L}$ is a partial differential operator with boundary operator $\mathcal{B}$, $w$ is the response of the system, $f, g$ are deterministic  and/or stochastic external forces and $\boldsymbol{\xi}(\omega), \omega\in\Omega\subset\mathbb{R}^{n_d}$ is a vector of $n_d$ random variables ($\omega$ is the random element on the probability space, which is omitted in the remaining of the paper). Since a GP  is a collection of random variables with a joint Gaussian distribution,  the continuous response $w $  can be interpreted as a realization of a Gaussian process  \cite{Krige51} over the parameter space $\boldsymbol{\xi}$
\begin{equation}\label{eq:Kriging}
w(\boldsymbol{\xi}) \sim \mathcal{GP}(\mu(\boldsymbol{\xi}), k(\boldsymbol{\xi}, \boldsymbol{\xi}^{'}))
\end{equation}
where $\mu(\boldsymbol{\xi})=\mathbb{E}[w(\boldsymbol{\xi}) ]$ and  $k(\boldsymbol{\xi}, \boldsymbol{\xi}^{'}) =\mathbb{E}[(w(\boldsymbol{\xi})-\mu(\boldsymbol{\xi}))(w(\boldsymbol{\xi}')-\mu(\boldsymbol{\xi}'))] $ are the mean and the covariance function, respectively. 


Suppose we have a training set composed of $N$ evaluations of our numerical model $\mathcal{M}$ denoted by $\mathcal{T}=\{(\boldsymbol{\Xi}, \vec{W})\}$ where  $\boldsymbol{\Xi}=\{\boldsymbol{\xi}^{(1)},\ldots,\boldsymbol{\xi}^{(N)}\}$ is the matrix of the input training parameters with corresponding response values  $\vec{W} = \{w(\boldsymbol{\xi}^{(i)}), \ldots,  w(\boldsymbol{\xi}^{(N)})\}$. Our objective is to predict the responses $\vec{W}^\star$ for a set of points $\boldsymbol{\Xi}^{\star}$ that are not contained in the training set. By definition, a  GP defines a prior over the function which can be converted into a posterior using the available  training data. If the training data have noise i.e. $y= w(\boldsymbol{\xi}) + \epsilon$ with $\epsilon \sim N(0,  \sigma_y^2)$ then, the GP model is not required to strictly interpolate the data, but should pass through the mean of observed data and provide a measure of the variance. Here we will not consider the case of noisy data, but will instead consider that our solution is precise at each training point. In this case, the prior joint density of the observed data and the test points is given by \cite{Rasmussen2006}
\begin{equation}
{{\vec{W}}\choose{\vec{W}^\star}} \sim N\bigg(\vec{0}, {{\vec{K} \quad \vec{K}_{\star}}\choose{\vec{K_{\star}^\intercal} \quad \vec{K}_{\star \star}}} \bigg) 
\end{equation}
where $\vec{K}=k(\boldsymbol{\Xi},\boldsymbol{\Xi})$, $\vec{K}_{\star}=k(\boldsymbol{\Xi}, \boldsymbol{\Xi}^\star)$ and $\vec{K}_{\star\star}=k(\boldsymbol{\Xi}^\star, \boldsymbol{\Xi}^\star)$. By conditioning the prior joint density on the observations $\vec{W}$, we obtain the posterior predictive density \cite{Rasmussen2006}:
\begin{eqnarray}
p(\vec{W}^\star|\boldsymbol{\Xi}^\star, \boldsymbol{\Xi}, \vec{W}) &\sim& N\left(\boldsymbol{\mu}_\star, \boldsymbol{\Sigma}_{\star}\right) \\
\boldsymbol{\mu}_\star &=&  \vec{K}_{\star}^\intercal\vec{K}^{-1}\vec{W} \\
\boldsymbol{\Sigma}_{\star} & = & \vec{K}_{\star\star} - \vec{K}_{\star}^\intercal\vec{K}^{-1}\vec{K}_{\star}
\end{eqnarray}
The above equations are easily extended for noisy data \cite{Rasmussen2006}. Hence, the approach proposed herein can be similarly extended for the case of simulations that generate solutions with additional embedded uncertainties.

The GP training model requires the selection of a suitable kernel for the coveriance function $k(\boldsymbol{\xi}, \boldsymbol{\xi}^{'})$. Numerous kernels are available, and the selection of the kernel is generally based on known (or assumed) properties of the response (e.g. continuity and smoothness of the function with respect to the input parameters). For our purposes, a Gaussian kernel is used to achieve a smooth, infinitely differentiable process. Additionally, one may also identify a trend associated with the mean value and apply regression techniques to solve for it's coefficients. See, e.g. \cite{schobi2015polynomial} who apply a polynomial chaos expansion for the trend. Regression of both the trend and the covariance function require estimation of the appropriate hyperparameters. Often, these are solved using a least squares technique (as is the case here) although other techniques can also be applied. The reader is referred to  \cite{Santner2003} for more details. 


\section{Data-driven GP regression on the Grassmann manifold}
\label{S:PM}

\begin{figure}[!ht]
	\centering
		\includegraphics[width=\textwidth]{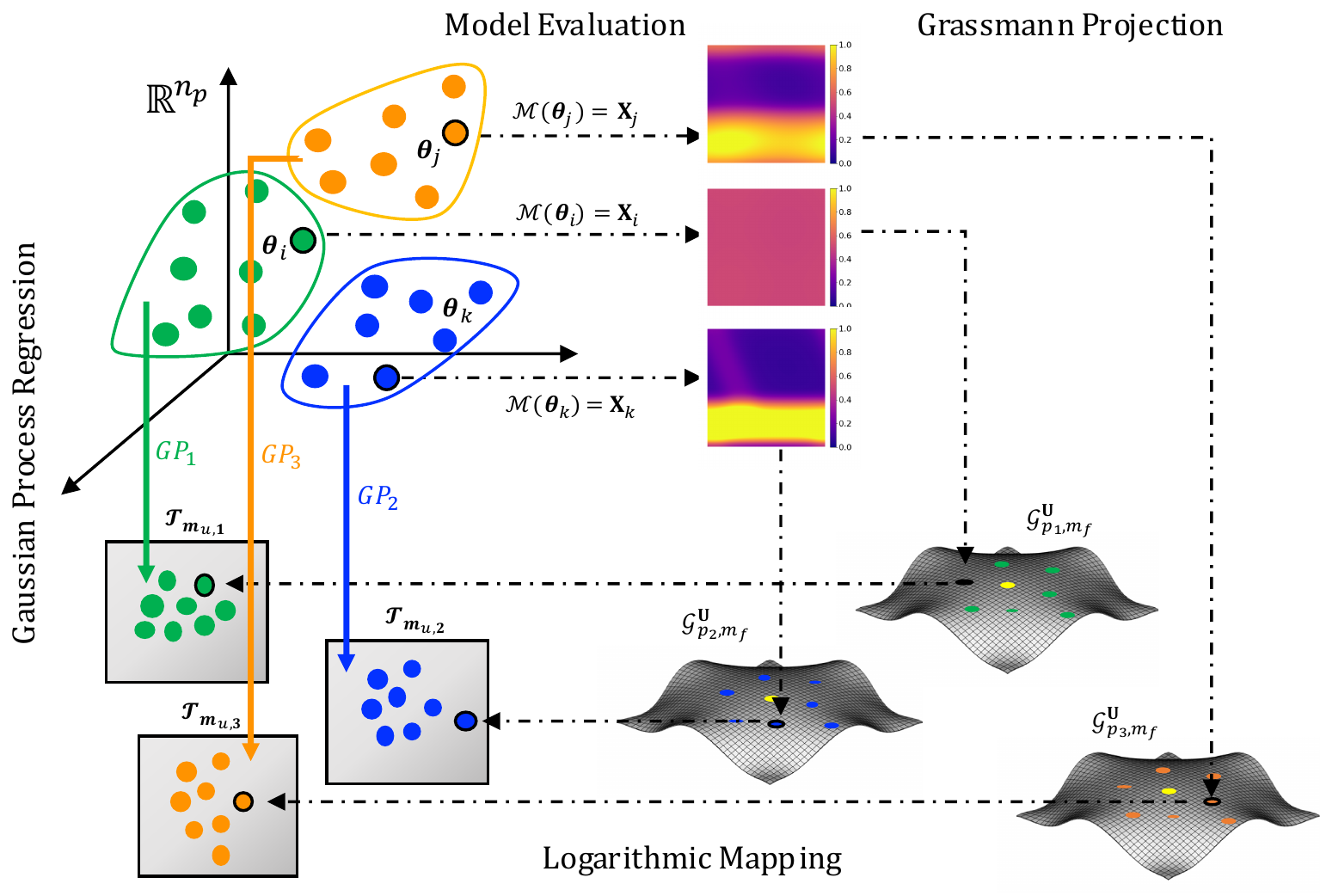}
	\caption{Flowchart of the proposed methodology.}
	\label{fig:framework}
\end{figure}

\noindent
In this section, we detail the proposed methodology for building GP surrogates for high-dimensional solutions on the Grassmann manifold. The process involves the following basic steps and is illustrated in Figure \ref{fig:framework}:
\begin{enumerate}
    \item \textbf{Training Simulations:} Run the training simulations and project the solutions from the training simulation onto the Grassmann manifold;
    \item \textbf{Solution Clustering:} Cluster the training simulations by similarity in their solutions;
    \item \textbf{Tangent Space Mapping:} Map each cluster onto the tangent space of the manifold defined at the Karcher mean of the cluster;
    \item \textbf{Parameter Space Clustering (Optional):} If necessary, perform a second clustering on the parameter space;
    \item \textbf{GP Regression \& Prediction:} Train a GP model for each cluster by defining the mapping from the parameter space (or a cluster in the parameter space) to the cluster tangent space. Employ the GP model for solution prediction at new points in the parameter space.
\end{enumerate}

\subsection{Training Simulations}
Consider a finite number $N_s$ of samples from the parameter space, i.e. realizations of the input random vector $\{\boldsymbol{\xi}_i\}_{i=1}^{N_s}$. These samples may be drawn using any arbitrary sampling method. For the purposes of illustration, we use simple Monte Carlo sampling although more advanced methods may improve the performance of the surrogate model. For each sample point, compute the high-dimensional (full-field) responses $\{w_i(\mathbf{x},t,\boldsymbol{\xi}_i)\}_{i=1}^{N_s}$ using the numerical model $\mathcal{M}$. The structure of the problem will depend on the nature of the problem. Dynamic problems having many degrees of freedom, solved using the finite element method for example, may take vector, matrix or tensor form. For consistency with the manifold projection algorithm, we recommend to recast the solution into matrix form $\{w_i(\mathbf{x},t,\boldsymbol{\xi}) \in \mathbb{R}^{n_{dof}}\}_{i=1}^{N_s}\rightarrow \{\vec{F}_i \in \mathbb{R}^{n_f\times m_f}\}_{i=1}^{N_s}$ where $n_{dof}=n_f \times m_f$. Vector-valued solutions can easily be reshaped as matrices, while tensor-valued responses can be recast using tensor unfolding methods. The selection of $n_f, m_f$ is decided by the user but a general rule-of-thumb used here is to make the matrix close to square, which often produces reduced solutions with the lowest total dimension. Note that matricization of the solution corresponds to a linear mapping that simply re-indexes the solution and therefore does not change the underlying dependency between components of the solution (e.g. degrees of freedom or time steps). For this reason, it is not expected to impact the approach significantly as long as the reshaping operation is performed consistently across all simulations. This is consistent with our observations below. 

Next, factorize each solution matrix $\{\vec{F}_i\}_{i=1}^{N_s}$ using a thin SVD as $\vec{F}_i=\vec{U}_i\boldsymbol{\Sigma}_i\vec{V}_i^\intercal$ where $\vec{U}_i \in \mathbb{R}^{n_f \times p_i}$, $\boldsymbol{\Sigma}_i \in \mathbb{R}^{p_i \times p_i}$ and $\vec{V}_i \in \mathbb{R}^{m_f \times p_i}$ and $p_i$ is the rank of $\vec{F}_i$. This factorization represents the linear projection of the solution onto the Grassmann manifold. Selection of the appropriate rank is based on standard truncation methods. In general, solutions corresponding to different input parameters may have different rank, i.e. $p_i \neq p_j$. This means that matrices (points) $\vec{U}_i$ and $\vec{U}_j$ may lie on different manifolds ($\mathcal{G}_{p_i, n_f}$ and $\mathcal{G}_{p_j, n_f}$). This will become important in the following step.

\subsection{Solution Clustering}
Given the local nature of the tangent space where interpolation is performed, it is critical that we attempted to build surrogates only over ``local'' regions of the Grassmann manifold where variations in the solution are considered small. Attempts to build surrogates across larger submanifolds will result in large errors associated with the logarithmic/exponential mapping and hence poor surrogate model predictions. In this section, we introduce a novel manifold clustering approach to identify ``local'' regions on the manifold.

The proposed approach is based on the spectral clustering technique, a widely-used approach for classification problems \cite{Luxburg2007}. One of the most widely used spectral clustering algorithms is the normalized cuts algorithm proposed in \cite{Shi2000}. Spectral clustering is performed by first constructing the so-called similarity (or affinity) matrix $\vec{W}=(w_{ij})$. In graph theory, one will also see this matrix referred to as the weight or adjacency matrix where the entries represent the weights of the graph edges. The similarity matrix $\vec{W}$ is a symmetric, positive definite matrix having strictly non-negative entries ($w_{ij}\ge0, \forall i,j$). Next, the degree of each point (vertex in the graph) is computed as $d_{ij} = \sum_{j=1}^{N_s} w_{ij}$ and degree matrix $\vec{D}=[d_{ij}]$ is constructed. From these matrices, the graph Lapacian is constructed. The graph Laplacian is the core component of the spectral clustering approach. Many variations of the graph Laplacians exist and have been extensively studies (as discusse in \cite{Luxburg2007}). Here, we use the symmetric normalized graph Laplacian defined as follows. First, the unnormalized graph Laplacian is defined by $\vec{L} = \vec{D}-\vec{W}$. The symmetric normalized graph Laplacian is given by \cite{Mohar1991,Mohar1997}:
\begin{equation}
    \vec{L}_{\text{sym}} = \vec{D}^{-1/2}\vec{L}\vec{D}^{-1/2} = \vec{I}- \vec{D}^{-1/2}\vec{W}\vec{D}^{-1/2}
\end{equation} 
Spectral clustering proceeds by computing the eigenvalues and eigenvectors of $\vec{L}_{\text{sym}}$. Defining $\boldsymbol{\Phi}$ as the matrix of eigenvectors, the rows of $\boldsymbol{\Phi}$ are then divided into a pre-specified number, $n_c$, of clusters using the k-means algorithm.

To apply spectral clustering on the Grassmann manifold, a few details of the algorithm are of particular importance. In particular, the weights of the similarity matrix, $w_{ij}$, must necessarily be produced by a positive semi-definite kernel possessing a valid distance metric. In conventional spectral clustering, it is sufficient to employ e.g. a Gaussian kernel based on a Euclidean distance. On the Grassmannian, this is a more challenging task. In Table \ref{tab:1}, we provide a brief list of distance measures. However, not all Grassmannian distance measures are valid metrics induced by a positive semi-definite kernel. Among the aforementioned distances, only the projection distance is a metric induced from a Grassmannian kernel \cite{Hamm:2008}. Note that other metrics exist, e.g.\ the Binet-Cauchy metric \cite{Hamm:2008}, but will not be discussed further here. The projection kernel, defined as:
\begin{equation}
    k_p(\vec{U}_i, \vec{U}_j) = ||\vec{U}_i^\intercal \vec{U}_j||_F^{2}
\end{equation}
is used herein to compute the similarity $w_{kl}$ and thus build the similarity matrix $\vec{W}$.

An additional important details is that the user must specify the number of clusters, $n_c$.  However, when comparing the solutions to high-dimensional nonlinear physics-based models, it is difficult to know a priori the  ``optimum'' number of solution clusters. This will be governed by the physical/mathematical nature of the problem. The solution may be governed by different regimes where different relevant mechanisms govern the response. For example, in a structural stability model, the displacement field may be governed by one of several buckling modes. We are not likely to know a priori what these modes are, or even how many modes play an important role.   

To determine the appropriate number of clusters, we propose an iterative procedure as a pre-processing step. Recall that our objective is to identify clusters whose solutions are all ``close'' on the manifold. Hence, in each iteration, we check that the solutions in each cluster are sufficiently close. If they are not, we increase the number of clusters. This process is repeated until our criterion is satisfied or we reach a user-defined maximum number of solution clusters, $n_{\max}$. We also require that each cluster contain a minimum number of points, $n_{min}$. A practical choice for $n_{\min}$ in order to ensure adequate surrogates over each cluster might be, e.g. $n_{\min}=10$, which establishes a bound of $n_{\max}=N_s/10$.

\begin{figure}[!htb]
 	\centering
 	\includegraphics[width=0.6\textwidth]{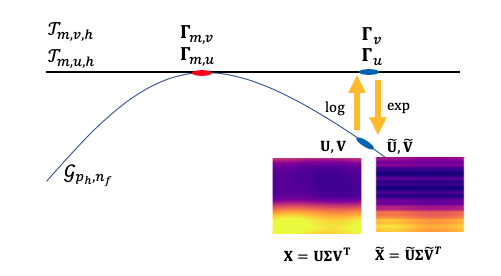}
 	\caption{Projection of a point from the manifold to the tangent space with origin the Karcher mean of the cluster and back onto the manifold. If the points are ``far'' away on the Grassmann the exponential mapping will introduce error in the solution.}
 	\label{fig:error}
 \end{figure}

To determine whether the points in each cluster are sufficiently close, we check the error introduced from the logarithmic mapping of each cluster as follows. For each solution cluster $C_h$, $h=1, \ldots, n_c$ having $N_h$ points, we  find the Karcher means $\mathbf{m}_{u,h}$ and $\mathbf{m}_{v,h}$ of points $ \{\vec{U}(\boldsymbol{\xi}^{(j)})\}_{j=1}^{N_h} \in C_h$  and $ \{\vec{V}(\boldsymbol{\xi}^{(j)})\}_{j=1}^{N_h} \in C_h$, respectively and project the points using the logarithmic mapping onto the corresponding tangent spaces (with origin at the Karcher mean). Then we project them back onto the Grassmann manifold using exponential mapping. This is expressed as follows:
\begin{subequations}
    \begin{equation}
    \{\vec{U}_j\in \mathcal{G}_{p_h, n_f}\}_{j=1}^{N_h} \rightarrow  \{\boldsymbol{\Gamma}_{u,j}\in \mathcal{T}_{\mathbf{m}_{u,h}}(\mathcal{G}_{p_h, n_f})\}_{j=1}^{N_h} \rightarrow 	\{\vec{\tilde{U}}_j\in \mathcal{G}_{p_h, n_f}\}_{j=1}^{N_h} 
    \end{equation}
    \begin{equation}
    \{\vec{V}_j\in \mathcal{G}_{p_h, m_f}\}_{i=1}^{N_h} \rightarrow  \{\boldsymbol{\Gamma}_{v,j}\in \mathcal{T}_{\mathbf{m}_{v, i}}(\mathcal{G}_{p_h, m_f})\}_{j=1}^{N_h} \rightarrow  	\{\vec{\tilde{V}}_j\in \mathcal{G}_{p_h, m_f}\}_{j=1}^{N_h}
    \end{equation}
\end{subequations}
where the indexes $u, v$ of $\boldsymbol{\Gamma}_{\cdot, j}$ and $\mathbf{m}_{\cdot,j}$ represent actions on the corresponding tangent spaces for the  left $\vec{U}_j$  and right eigenvector $\vec{V}_j$,  respectively. In these equations $p_h=\max( p_j:\{\vec{U}_j\in \mathcal{G}_{p_j, n_f}\}_{j=1}^{N_h})$ is the maximum rank of the points in the cluster. If the points are ``far'' away on the Grassmannian, then the projection from the manifold to the tangent space and back will induce significant error in the solution (see Fig.\ \ref{fig:error}). Numerically the point-wise projection error is quantified with the Frobenius norm
\begin{equation}
\alpha_j = || \vec{U}_j\boldsymbol{\Sigma}_j\vec{V}_j^\intercal - \vec{\tilde{U}}_j\boldsymbol{\Sigma}_j\vec{\tilde{V}}_j^\intercal ||_F \\
\end{equation}
The average projection error of each cluster is then estimated by
\begin{equation}\label{error_mse}
  \epsilon^h = \frac{1}{N_h}\sum_{j=1}^{N_h}  \alpha_i
 \end{equation} 
If a large fraction of the clusters have an average projection error less than a threshold, e.g.\ $\epsilon^h\leq 10^{-3}$, then the iterations stop and we proceed with $n_c$ clusters. Otherwise we increase the number of clusters by $n_c = n_c + 1$ and repeat the procedure. A typical fraction value is taken as $0.8-0.9$, meaning that we allow 10--20\% of the clusters to introduce some error given the difficulty of achieving small error across the full parameter space from a limited number of simulations. Of course, the user may choose to be more strict with the error criterion and/or clusters that produce significant errors can be flagged as erroneous or marked for potential improvement.

\subsection{Tangent Space Mapping}
Once the clusters are identified, it is necessary to map the right and left basis vectors of each solution in each cluster onto the tangent space of that cluster ($\vec{U}_j$ and $\vec{V}_j$ lie on different manifolds) defined at the cluster Karcher mean, $\mathcal{T}_{\mathbf{m}_h}$. This is performed using the logarithmic mapping given in Eq.\ \eqref{eq:8}. That is, we perform the following operations
\begin{subequations}
    \begin{equation}
    \boldsymbol{\Gamma}_{u,j} = \log_{\mathcal{T}_{\mathbf{m}_{u,h}}}(\vec{U}_j)
    \end{equation}
    \begin{equation}
    \boldsymbol{\Gamma}_{v,j} = \log_{\mathcal{T}_{\mathbf{m}_{v,h}}}(\vec{V}_j)
    \end{equation}
\end{subequations}
where, again $\vec{F}_j=\vec{U}_j\boldsymbol{\Sigma}_j\vec{V}_j^\intercal$. This step maps the solutions to points in a space where interpolation can be performed and therefore a surrogate model can be constructed.

\subsection{Cluster Gaussian Process Regression}
For each cluster, the surrogate modeling scheme requires the construction of three GPs. This is necessary because the SVD of a given solution provides two sets of basis vectors (the left and right singular vectors, $\vec{U}$ and $\vec{V}$ respectively) and a set of scale factors (the singular values, $\sigma$). To reconstruct the full solution at a new point in the parameter space, each of these components must be predicted. We therefore construct the following GP models for the tangent space projections $\boldsymbol{\Gamma}_{u,j}$, $\boldsymbol{\Gamma}_{v,j}$ and the singular values $\boldsymbol{\Sigma}_j$ for each cluster
\begin{subequations}
 	\begin{equation}
        \tilde{\boldsymbol{\Gamma}}_{u,j}(\boldsymbol{\xi}) \sim \mathcal{GP}(\mu_{u,j}(\boldsymbol{\xi}), k_{u,j}(\boldsymbol{\xi}, \boldsymbol{\xi}^{'}))
 	\end{equation}
 	\begin{equation}
        \tilde{\boldsymbol{\Gamma}}_{v,j}(\boldsymbol{\xi}) \sim \mathcal{GP}(\mu_{ v, j}(\boldsymbol{\xi}), k_{v, j}(\boldsymbol{\xi}, \boldsymbol{\xi}^{'}))
 	\end{equation}
 	\begin{equation}
 	 	 \tilde{\boldsymbol{\Sigma}}_j(\boldsymbol{\xi}) \sim \mathcal{GP}(\mu_{\boldsymbol{\sigma}_j}(\boldsymbol{\xi}), k_{\boldsymbol{\sigma}_j}(\boldsymbol{\xi}, \boldsymbol{\xi}^{'}))
 	\end{equation}
 \end{subequations}
Here, prediction of these matrix quantities occurs component-wise. Given the structure of the data, one could envision improving these predictors with multi-variate Kriging models / coKriging although we find that sufficient accuracy can be achieved without this.

By training these local GPs we can predict the full response field at a new point $\boldsymbol{\xi}^\star$ in the parameter space very efficiently. To do so, we must first identify which solution cluster $\boldsymbol{\xi}^\star$ belongs to. This is done by simply finding the nearest point in the parameters space and assigning the new point to the same cluster as its neighbor. To identify the nearest point, we use a Euclidean distance. However, one may consider other distance measures based on different metrics. Additionally, distances for input possessing non-uniform distribution may require an isoprobabilistic transformation prior to distance computation. The approximate solution is then constructed from the surrogate predictions by first performing an exponential mapping of the estimated tangent space representations:
\begin{subequations}
    \begin{equation}
    \tilde{\vec{U}}_j = \exp_{\mathcal{T}_{\mathbf{m}_{u,h}}}(\tilde{\boldsymbol{\Gamma}}_{u,j})
    \end{equation}
    \begin{equation}
    \tilde{\vec{V}}_j = \log_{\mathcal{T}_{\mathbf{m}_{v,h}}}(\tilde{\boldsymbol{\Gamma}}_{v,j})
    \end{equation}
\end{subequations}
and multiplying the matrices to assemble to prediction as
\begin{equation}
\tilde{\vec{F}}(\boldsymbol{\xi}^\star) = \tilde{\vec{U}}\tilde{\boldsymbol{\Sigma}}\tilde{\vec{V}}^\intercal.
\end{equation}






\subsection{Optional Parameter Space Clustering}

\noindent
In some cases, clustering the points of the training set according to solution similarity may result in clusters whose points lie in disjoint regions areas of the parameter space. This is often a byproduct of significant nonlinearity in the model and is related to the challenge of non-uniqueness in model inversion. Practically, this means that points that are ``far'' apart in the parameter space map to the same (or very similar) solutions, they become clustered together, and the points between them lie in another cluster. This becomes problematic only when the points in the parameter space lie in disjoint regions because continuity in the GP model is lost.  In many such cases, training a GP surrogate model to such a data set will result in a poor approximation.


In such cases, we propose a second-level of clustering in the parameter space. That is, within each solution cluster we perform an additional ``input subclustering.'' This subclustering follows a more conventional approach as the parameter space is generally a Euclidean space having the usual metrics of distance. In our implementation, we utilize the density-based spatial clustering of applications with noise (DBSCAN) algorithm proposed in \cite{DBSCAN1996}, with the objective of grouping the cluster training points based on their location in the parameter space.  Then, for each of the input subclusters we can train a GP using the described procedure.

\section{Numerical examples}
\label{S:Examples}

In this section, we apply the proposed approach for three examples that are intended to explore its performance and limitations. In the first example, we apply it to a model system of nonlinear ordinary differential equations (ODEs). We develop a surrogate model for the time history response of one component of the system. The second and third examples consider a model of plastic deformation in amorphous materials. One example considers a low-dimensional input parameter space that governs the mean and variance of the stochastic initial conditions to enable visualization of the proposed method. The other example considers the problem in a multi-scale setting where the input parameter vector is larger.

\subsection{Example 1: Kraichnan-Orszag SODE}
	
The Kraichnan-Orszag (KO) three mode problem is a set of nonlinear three-dimensional stochastic ordinary differential equations (SODEs) previously investigated by several authors \cite{KO1967, MEgPC2005,MEgCM2006,MEPCM2008,ASGC2009}. In accordance with \cite{MEgPC2005}, we rotate the equations by $\pi/4$ around the $w_{3}$-axis in the phase space, such that the system is defined as follows:
\begin{eqnarray}\label{KOR}
\frac{dv_{1}}{dt}&=&v_{1}v_{3} \nonumber \\
\frac{dv_{2}}{dt}&=&-v_{2}v_{3} \\
\frac{dv_{3}}{dt}&=&-v_{1}^{2}+v_{2}^{2} \nonumber 
\end{eqnarray}
The equations are subject to random initial conditions $v_{1}(0) = v_{1}(0,\omega), \quad v_{2}(0) = v_{2}(0,\omega), \quad v_{3}(0) = v_{3}(0,\omega)$, where $\omega$ indexes the probability space. For the purpose of illustration, we select the initial condition $v_1(0)$ to be deterministic and equal to 1.0 while the initial conditions $v_2(0)$ and $v_3(0)$ are considered random and are defined as $v_{2}(0) = 0.1\times\xi_{1}$ and $v_{3}(0) =\xi_{2}$, respectively where $\xi_{1}, \xi_{2}$ are uniformly distributed random variables over the range $(-1, 1)$. The solution to the system is periodic, but having discontinuity at $v_{1}=0$ corresponding to a solution with infinite period.

\begin{figure}[!htb]
	\centering
	\includegraphics[width=0.7\textwidth]{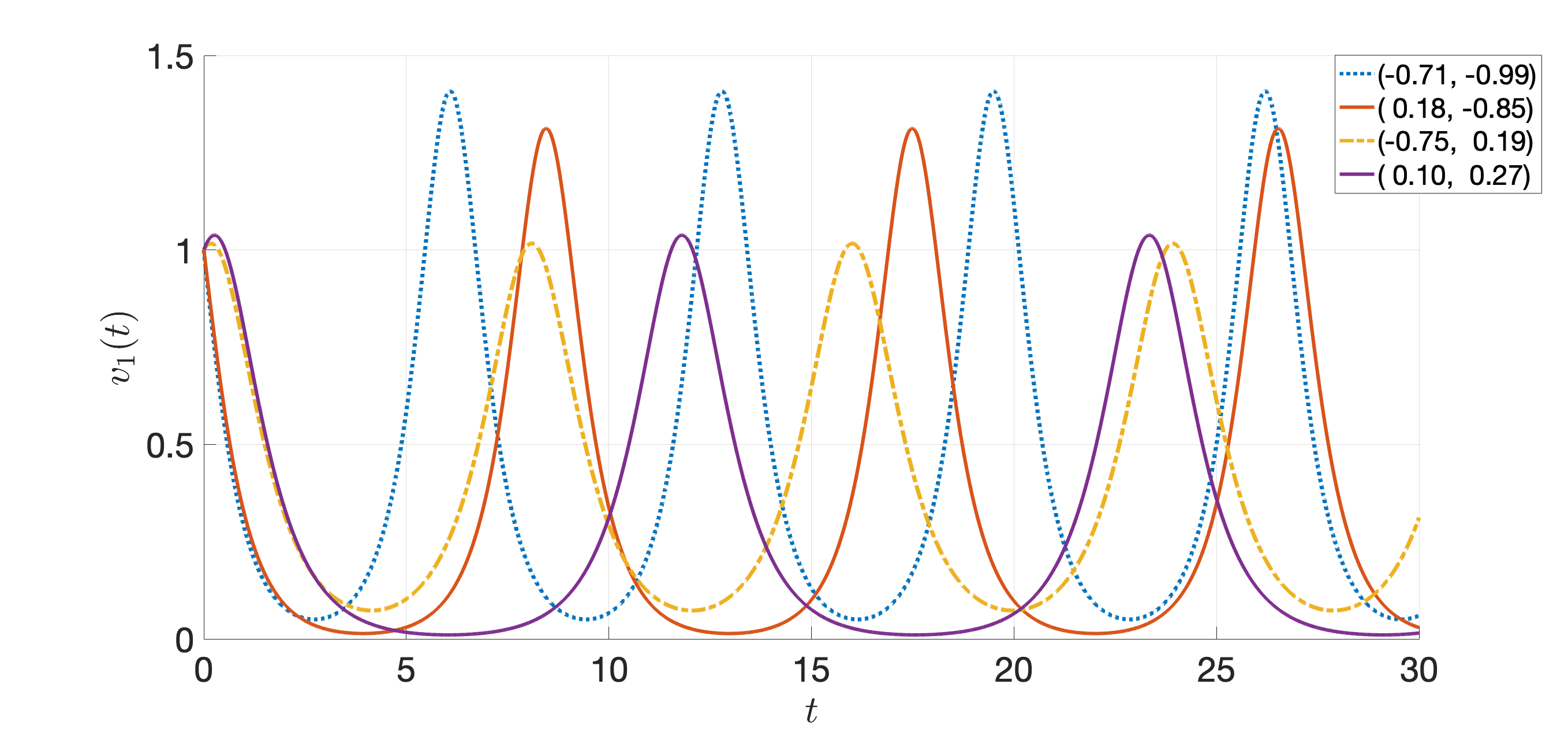}
	\caption{The solution $v_1(\vec{t})$ for the KO problem with two random variables. Different realizations of the initial conditions $v_2(0)$ and $v_3(0)$ lead to different behavior in the solution.}
	\label{fig:gfun_ko2d}
\end{figure}

The SODE system is solved for a period of T = 30 with a time step $\Delta t =  0.003$ which results in a response vector $w(\vec{t}, \boldsymbol{\xi})= [v_1(\vec{t},  \boldsymbol{\xi})), v_2(\vec{t},  \boldsymbol{\xi})), v_3(\vec{t},  \boldsymbol{\xi}))] \in \mathbb{R}^{3\times10^4}$.   For the purpose of this work we generated $N_s=1024$ realizations of the parameter vector $\boldsymbol{\xi}=(\xi_1, \xi_2)$ using Monte Carlo simulation and solve the system for the corresponding initial conditions. To illustrate the performance of the surrogate for this problem we study only the first component of the solution $w(\vec{t}, \boldsymbol{\xi})=v_1(\vec{t},  \boldsymbol{\xi})) \in \mathbb{R}^{10^4}$. The method can easily be applied to all three components, but we refrain from doing so in order to more directly illustrate the surrogate. Figure \ref{fig:gfun_ko2d} shows the sensitivity of the solution to variations of the initial conditions $v_2(0)$ and $v_3(0)$, where four different realizations of the random vector $\boldsymbol{\xi}=(\xi_1, \xi_2) =(-0.71, -0.99)$, $(0.18, -0.85)$, $(-0.75, 0.19)$, and $(0.10, 0.27)$ result in different behavior in the response.

\begin{figure}[!htb]
	\centering
	\begin{subfigure}[t]{0.25\textwidth}
		\centering
		\includegraphics[width=\textwidth]{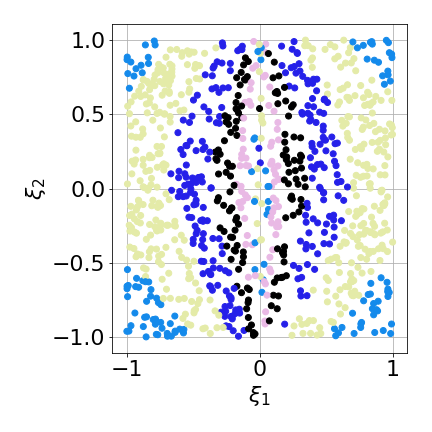}
		\caption{}
	\end{subfigure}%
	\begin{subfigure}[t]{0.25\textwidth}
		\centering
		\includegraphics[width=\textwidth]{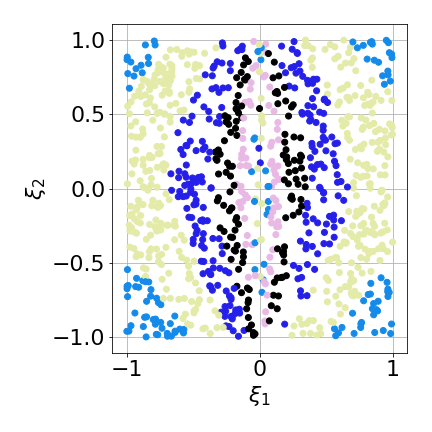}
		\caption{}
	\end{subfigure}%
	\begin{subfigure}[t]{0.25\textwidth}
	    \centering
	    \includegraphics[width=\textwidth]{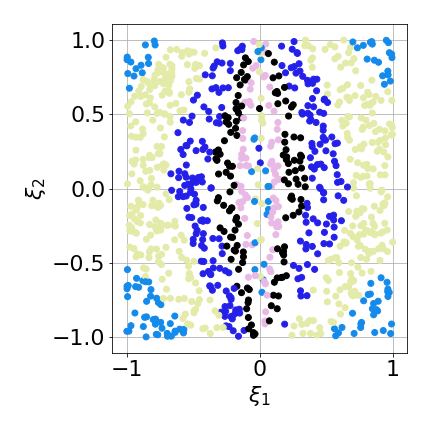}
	    \caption{}
    \end{subfigure}%
    \begin{subfigure}[t]{0.25\textwidth}
	    \centering
	    \includegraphics[width=\textwidth]{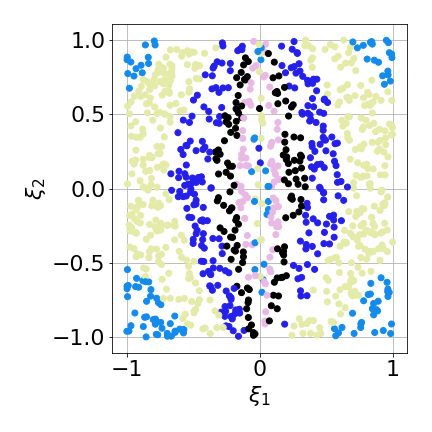}
	    \caption{}
    \end{subfigure}%
	\caption{Solution clustering patterns of  $N_s=1000$ input vectors $\boldsymbol{\xi}$ for  $n_c=5$ clusters and size of matrix (a) $\vec{F}_{1, i}\in \mathbb{R}^{100 \times 100}$,  (b) $\vec{F}_{1, i}\in \mathbb{R}^{200 \times 50}$,  (c) $\vec{F}_{1, i}\in \mathbb{R}^{1000 \times 10}$ and (d) $\vec{F}_{1, i}\in \mathbb{R}^{2000 \times 5}$.}
	\label{fig:comp_cl}
\end{figure}

Initially, one needs to cast the vector-valued time history response of each solution into a matrix form, denoted $\vec{F}_i$.  To highlight the insensitivity of the method to  matrix shape, consider four different cases,  $v_{1, i} \rightarrow \vec{F}_{i}\in \mathbb{R}^{100 \times 100}(n_f=100, m_f=100)$,  $v_{1, i} \rightarrow \vec{F}_{i}\in \mathbb{R}^{200 \times 50}(n_f=200, m_f=50)$, $v_{1, i} \rightarrow \vec{F}_{i}\in \mathbb{R}^{1000 \times 10}(n_f=1000, m_f=10)$ and $\vec{F}_{1, i}\in \mathbb{R}^{2000 \times 5}(n_f=2000, m_f=5)$. Using these reshaped solutions, we perform solution clustering with a small number of clusters $n_c=5$.  Figure \ref{fig:comp_cl} shows that the solution clustering is independent of the matrix size. For the remainder of this example, $\vec{F}_{i}\in \mathbb{R}^{100 \times 100}$ is selected.

\begin{figure}[!htb]
	\centering
	\begin{subfigure}[t]{0.25\textwidth}
		\centering
		\includegraphics[width=\textwidth]{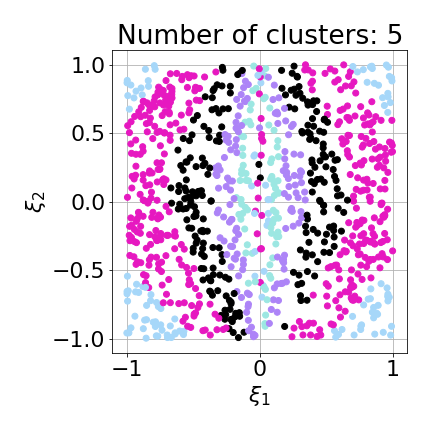}
		\caption{}
	\end{subfigure}%
	\begin{subfigure}[t]{0.25\textwidth}
		\centering
		\includegraphics[width=\textwidth]{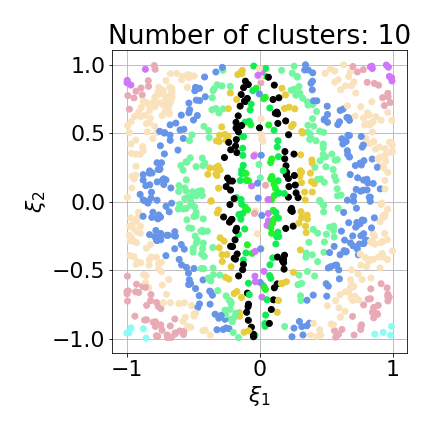}
		\caption{}
	\end{subfigure}%
	\begin{subfigure}[t]{0.25\textwidth}
	    \centering
	    \includegraphics[width=\textwidth]{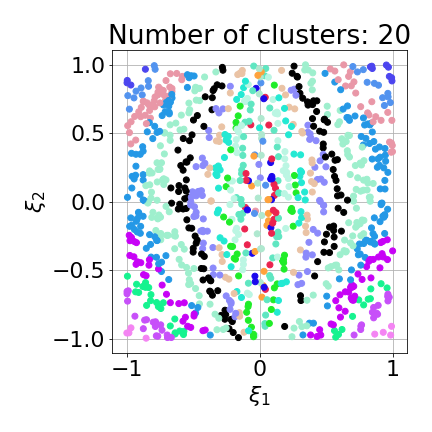}
	    \caption{}
    \end{subfigure}%
    \begin{subfigure}[t]{0.25\textwidth}
	    \centering
	    \includegraphics[width=\textwidth]{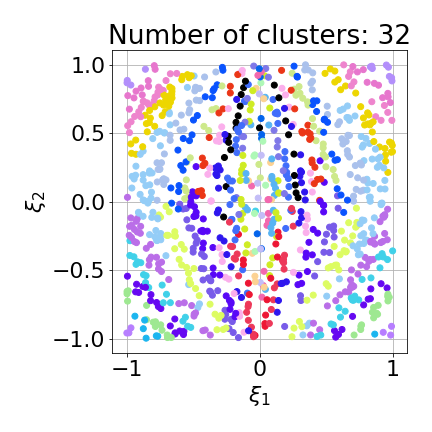}
	    \caption{}
    \end{subfigure}%
	\caption{Clustering patterns of the $N_s=1000$  MCS samples of $\boldsymbol{\xi}$ using (a) $n_c=5$, (b) $n_c=10$, (c) $n_c=20$ and (d) $n_c=32$ clusters.}
	\label{fig:clusters}
\end{figure}

\begin{figure}[!htb]
	\centering
	\begin{subfigure}[t]{0.33\textwidth}
		\centering
		\includegraphics[width=\textwidth]{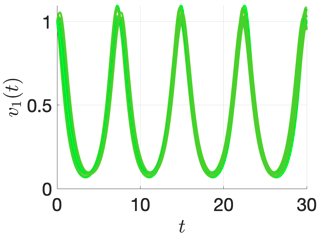}
		\caption{}
	\end{subfigure}%
	\begin{subfigure}[t]{0.33\textwidth}
		\centering
		\includegraphics[width=\textwidth]{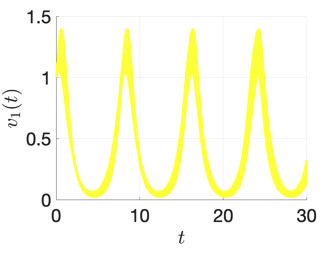}
		\caption{}
	\end{subfigure}%
	\begin{subfigure}[t]{0.33\textwidth}
	    \centering
	    \includegraphics[width=\textwidth]{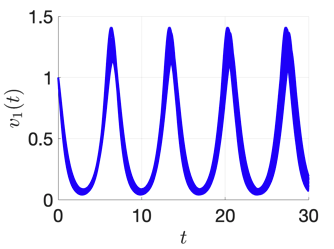}
	    \caption{}
    \end{subfigure}%
	\caption{High-dimensional $v_1(\vec{t})$ solutions that belong to three from a total of $n_c=32$ clusters.}
	\label{fig:gfun_cluster}
\end{figure}

A key component of the proposed method is the solution clustering built from distance metrics induced from appropriate Grassmannian kernels. Figure \ref{fig:clusters} illustrates the solution clustering for the Monte Carlo training points in the parameter space with different numbers of clusters, $n_c=5, 10, 20$ and $32$. We observe that there exist symmetrical bands in the parameter space whose time history response is similar. In Figure \ref{fig:gfun_cluster}, the corresponding solutions are shown for three different clusters from a total of 32 clusters.

The efficiency of the proposed method is highly dependent on the number of selected solution clusters, $n_c$. If small variations in the parameter space cause sharp changes in the solution then a large number of clusters may be required. However, knowledge of the number of clusters is required a priori. Following the proposed procedure, we estimate the ``optimal'' number of clusters is $n_c = 32$ under the requirement that $95\%$ of the clusters have a threshold value less than or equal to $10^{-3}$ in Eq.\ \eqref{error_mse}. 



For each of the 32 clusters, we need to perform a second subclustering in the parameter space because, as we can see from Figure \ref{fig:clusters}, the symmetry plane causes clusters to be disjointly connected to their reflection over the $\xi_1=0$ plane. Additional disjoint clusters can also be seen in the $\xi_2$ direction for a small number of clusters.

Next, for each subcluster we construct the Gaussian process models. As previously mentioned, we selected the Gaussian kernel, also known as the a radial basis function (RBF) kernel or squared exponential kernel, for the covariance. This kernel is define as
\begin{equation}
k(x_i, x_j) =\exp\left(-\frac{(x_i - x_j)^2}{2l^2}\right)
\end{equation}
\noindent
which is parameterized by a scale-parameter\footnote{The scale parameter $l$ can either be a scalar (isotropic variant of the kernel) or a vector with the same number of dimensions as the inputs (anisotropic variant of the kernel).}  $l>0$.  The Gaussian kernel is infinitely differentiable, which implies that GPs with this kernel have mean square derivatives of all orders, and are thus very smooth. The hyperparameter $l$ is cluster dependent and solved through a least squares optimization with initial value $l=1.0$.  

\begin{figure}[!htb]
	\centering
	\includegraphics[width=1.0\textwidth]{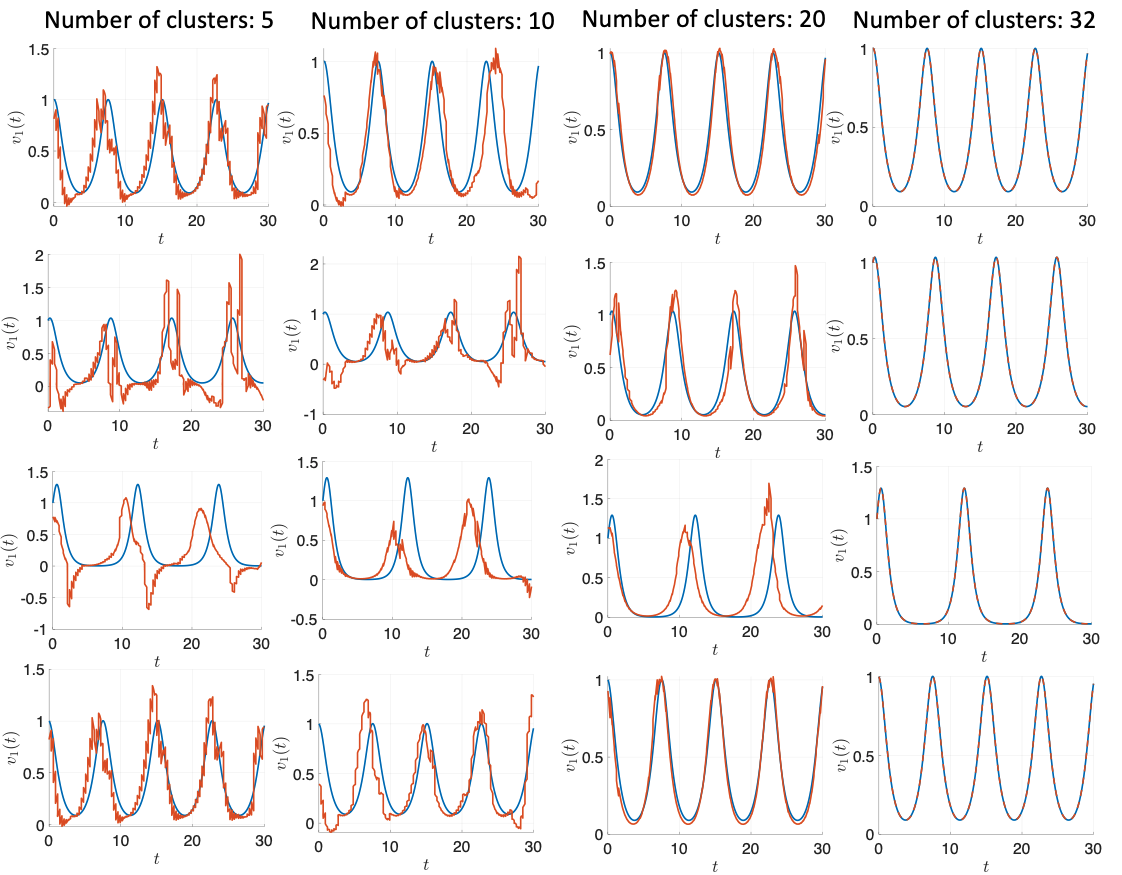}
	\caption{Exact (blue line) vs approximate (red line) solution of $v_1(\vec{t}, \boldsymbol{\xi})$ for different number of solution clusters. From top to bottom each row corresponds to point $\boldsymbol{\xi}_1=(0.89, -0.07)$, $\boldsymbol{\xi}_2=(-0.03, 0.81)$,  $\boldsymbol{\xi}_3=(-0.53, 0.26)$  and $\boldsymbol{\xi}_4=(0.91, 0.02)$.}
	\label{fig:point1}
\end{figure}

In order to visualize the quality of the approximate solution obtained from the Gaussian processes we arbitrarily select four points ($\{\boldsymbol{\xi}^\star\}$ = $(0.89, -0.07)$, $(-0.03, 0.81)$, $(-0.53, 0.26)$ and  $(0.91, 0.02)$) in the parameter space and compare the exact solution with the GP approximation for increasing number of clusters $n_{c} =5, 10, 20$ and the optimal $n_c=32$.  We observe from Figure  \ref{fig:point1} that the surrogate model performance improves significantly as we increase the number of clusters up to the optimal number. This is because, for small number of clusters, significant error is introduced from the logarithmic mapping when points are far apart on the manifold.

\begin{figure}[!htb]
	\centering
	\includegraphics[width=0.30\textwidth]{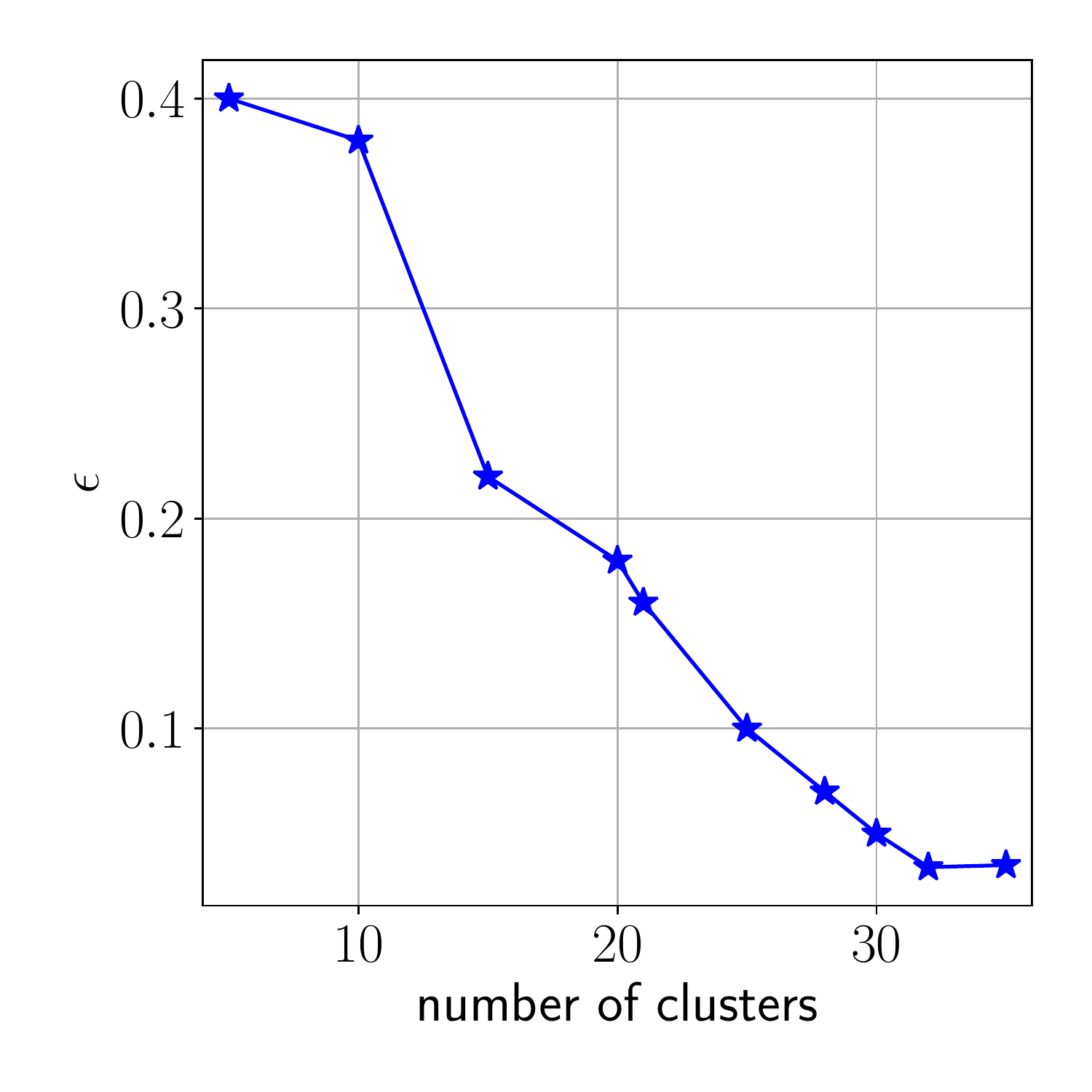}
	\caption{Surrogate model solution error as a function of the number of clusters $n_c$.}
	\label{fig:rmse}
\end{figure}

\begin{figure}[!htb]
	\centering
	\begin{subfigure}[t]{0.25\textwidth}
		\centering
		\includegraphics[width=\textwidth]{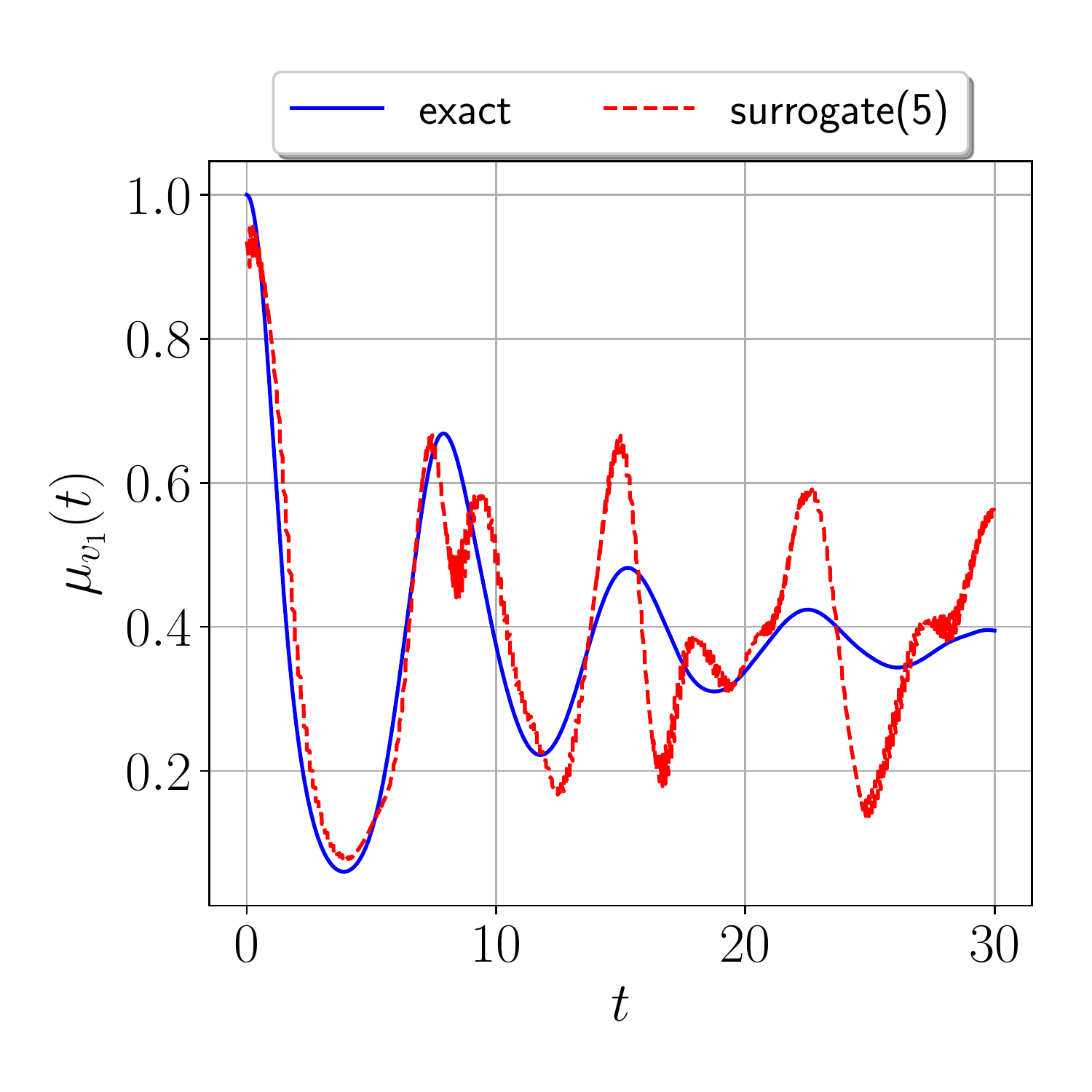}
		\caption{}
	\end{subfigure}%
	\begin{subfigure}[t]{0.25\textwidth}
		\centering
		\includegraphics[width=\textwidth]{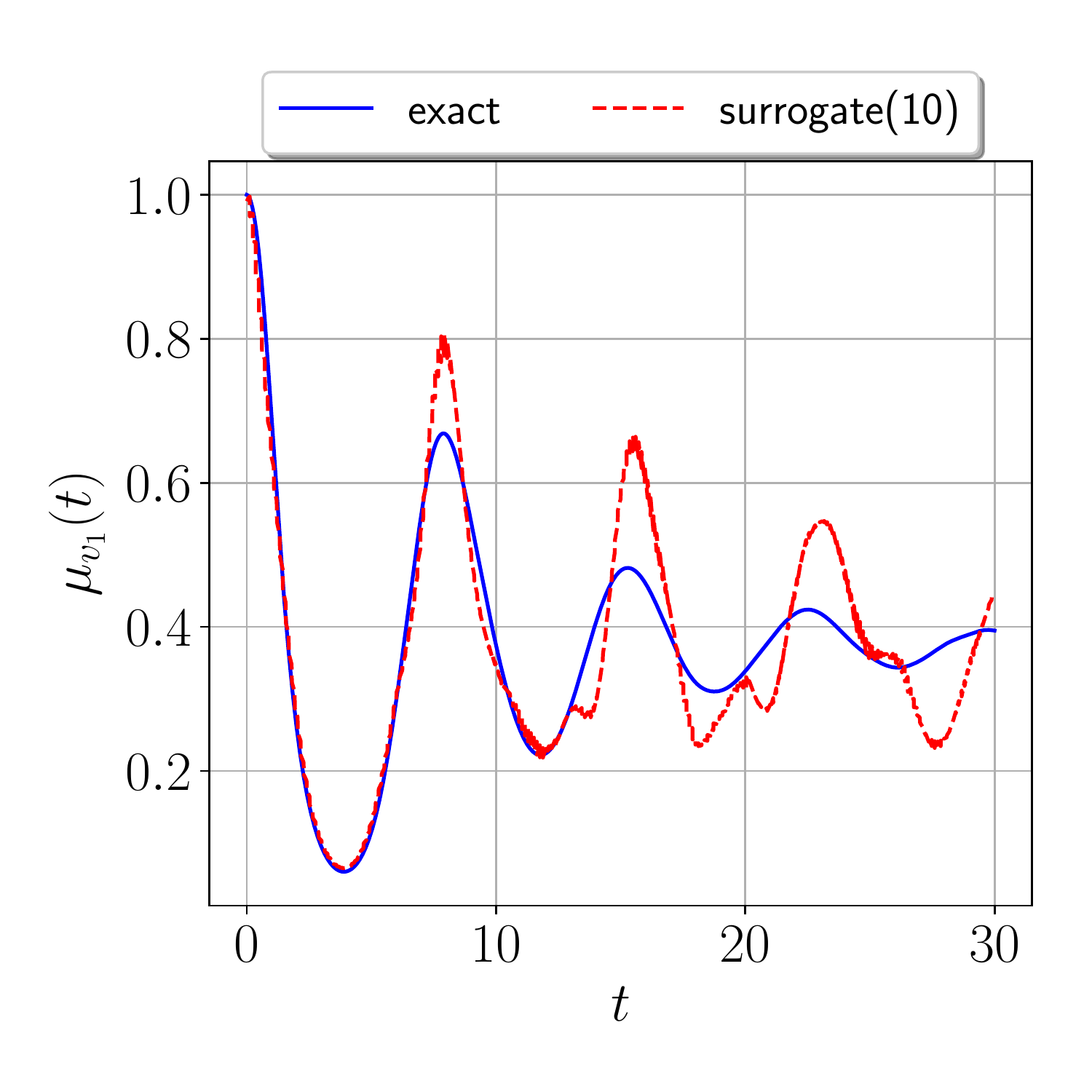}
		\caption{}
	\end{subfigure}%
    \begin{subfigure}[t]{0.25\textwidth}
	    \centering
	    \includegraphics[width=\textwidth]{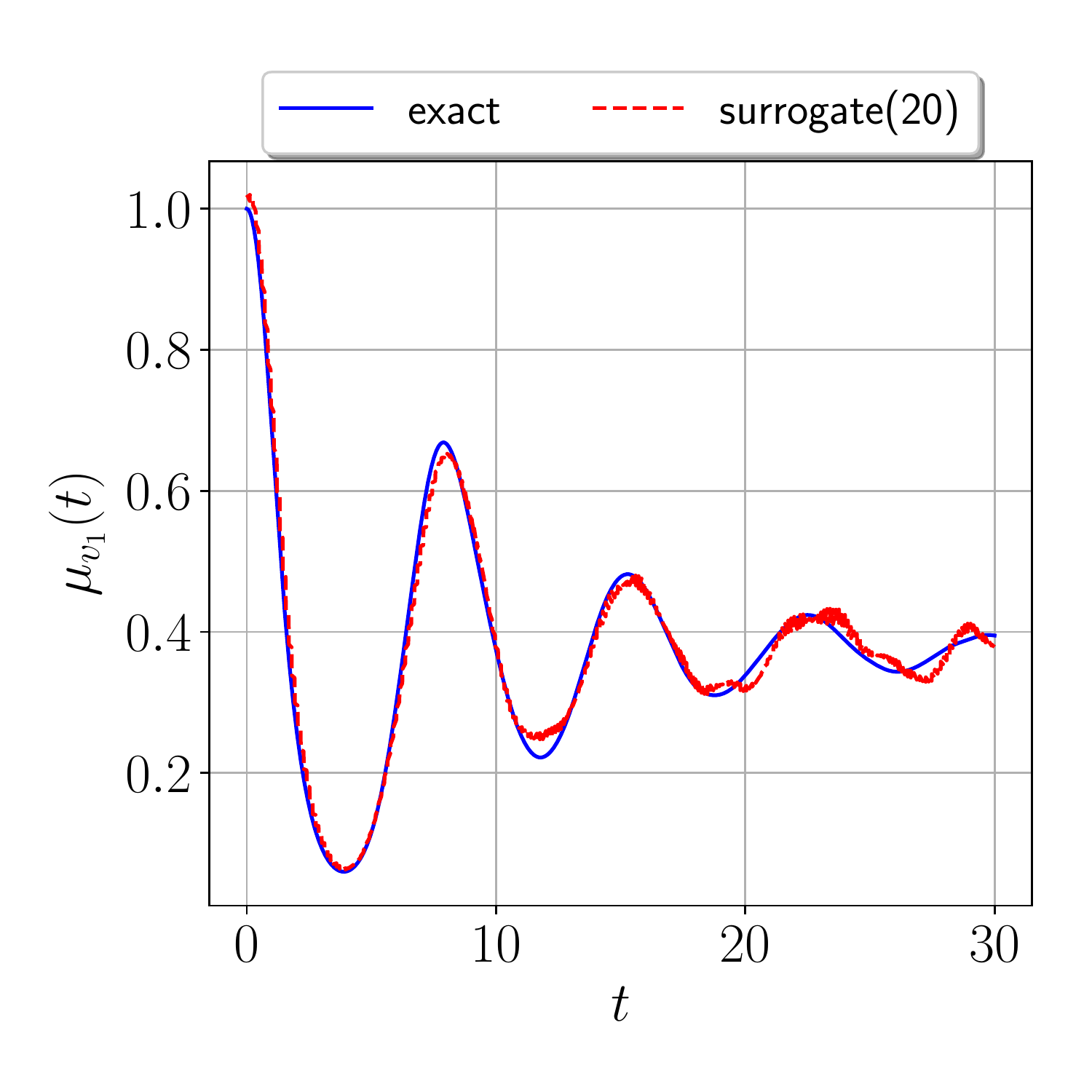}
	    \caption{}
    \end{subfigure}%
    \begin{subfigure}[t]{0.25\textwidth}
	    \centering
	    \includegraphics[width=\textwidth]{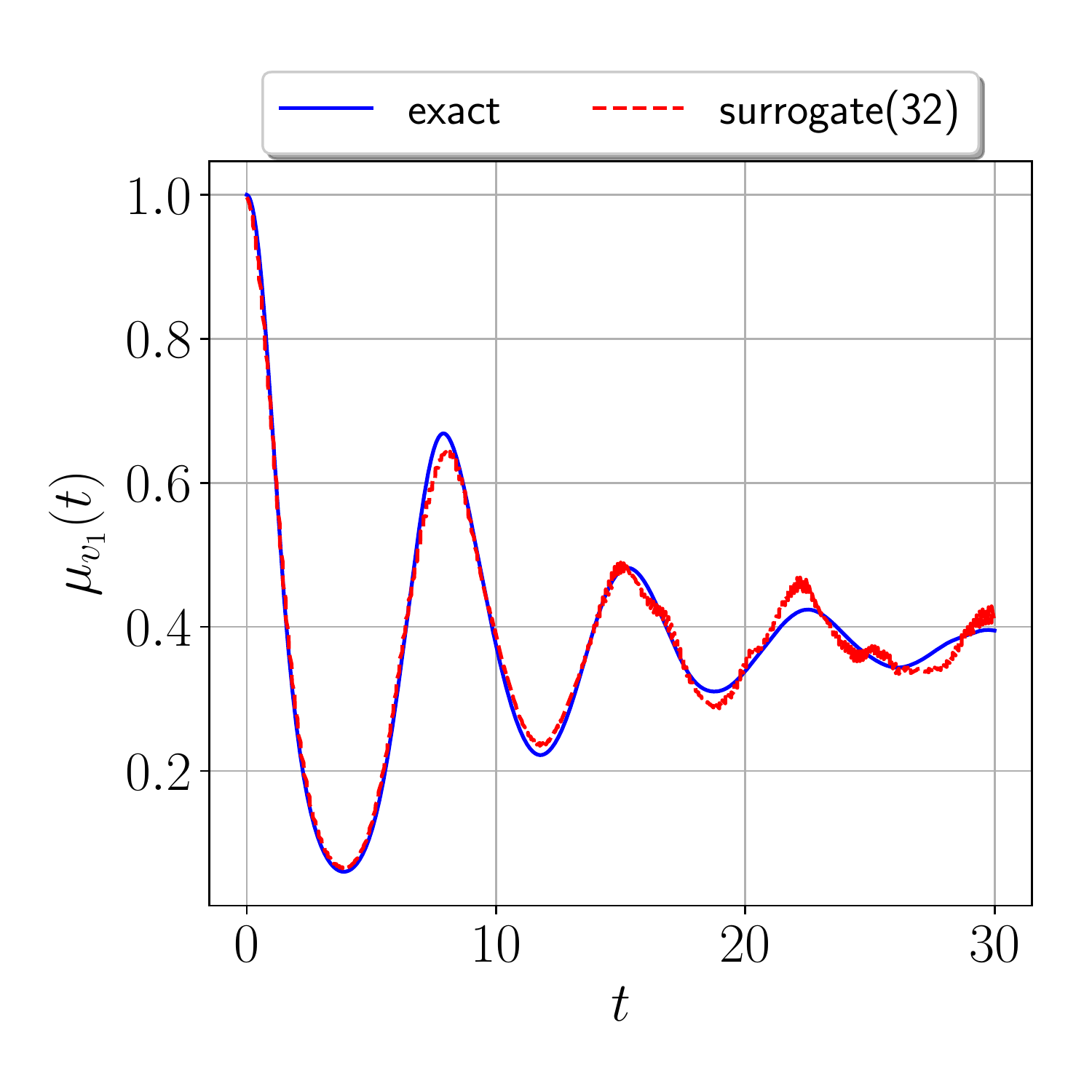}
	    \caption{}
    \end{subfigure}%
	\caption{ Mean of ${v}_1$ (red dashed line) vs. time for (a) $n_{cl}=5$, (b) $n_{cl}=10$, (c) $n_{cl}=20$  and (d) $n_{cl}=32$, compared to the exact mean $\mu_{v_1}(t)$ (solid blue line).}
	\label{fig:mean_w}
\end{figure}

To quantify the quality of the surrogate results, we generated  $N_{\text{test}}=3000$ additional parameter realizations and calculated the error between the true solution and the approximate solution from the GP using the following error metric
\begin{equation}\label{error_metric}
\epsilon = \frac{1}{N_{\text{test}}}\sum_{i = 1}^{N_{\text{test}}} ||A_i||_F, \quad \text{for}~i = 1, \ldots, 3000
\end{equation}
where  $||\cdot ||_F$ is the  Frobenius norm of the matrix  $\vec{A}_i = \vec{F}_{i} - \tilde{\vec{F}}_{i}$. Figure \ref{fig:rmse} plots this error as a function of the number of clusters.  Additionally, Figure \ref{fig:mean_w} shows the mean $\mu_{v_1}(t)$ from the $N_s=3000$ test points calculated using the proposed local surrogate solution and the true solution. Again, ss the number of clusters increases, the surrogate performance improves such that we obtain a near exact match for the optimal clustering.

\subsection{Example 2: Continuum modeling of plasticity in an amorphous solid}

\noindent
The proposed surrogate modeling approach is demonstrated here on a continuum hypo-elastoplastic material model for amorphous solids based on the shear transformation zone (STZ) theory of plasticity \cite{Falk1998, Langer2008, Falk2011,Bouchbinder2009, Bouchbinder2009b, Bouchbinder2009c}. The model is solved using an Eulerian finite difference numerical integration that leverages a quasi-static approximation to allow for very large deformations that cannot be realized in Lagrangian schemes that are typical for modeling solids \cite{Rycroft2015, boffi2020parallel}.  

\subsubsection{STZ Theory of Plasticity}
In amorphous solids, an emerging theory is that irreversible plastic deformation is mediated through atomic rearrangements in small clusters of atoms referred to as shear transformation zones (STZs) \cite{Falk1998}. When large shear stresses are applied to the material, STZs rearrange in localized regions that coalesce into thin bands of large local deformation called shear bands. The STZ theory of plasticity aims to connect these local rearrangements to larger-scale plastic deformation through an effective temperature that correlates with the density of STZs. The effective temperature is a measure of structural disorder defined as 
\begin{equation}\label{eq:31}
T_{eff} =\dfrac{\partial U_c}{\partial S_c}
\end{equation}
where $U_c$ and $S_c$ are the potential energy and the entropy of the configurational degrees of freedom under the assumption that total energy and total entropy can be separated into kinetic and configurational components as $U=U_c+U_k$ and $S=S_c+S_k$, respectively. A nondimensionalized effective temperature can then be defined as
\begin{equation}\label{eq:31}
\chi = \frac{k_BT_{eff}}{E_Z}
\end{equation}
$k_B$ is the Boltzmann factor and $E_Z$ is the STZ formation energy.

The effective temperature is a spatially varying quantity in a given material specimen and the STZ theory of plasticity defines a set of equations that describe its evolution along with the evolution of plastic deformation caused by an applied stress. The theory involves two essential equations. The first is a plastic flow rule that relates the plastic rate of deformation tensor to the effective temperature and the stress. Several versions of this have been proposed, one of which is given by:
\begin{equation}\label{eq:32}
    \mathbf{D}^{pl}=\dfrac{1}{\tau_0}\exp\left\{-\left(\dfrac{-e_z}{k_B\chi}+\dfrac{\Delta}{k_B T}\right)\right\}\cosh\left(\dfrac{\Omega\epsilon_0\bar{\sigma}}{k_BT}\right)\left(1-\dfrac{\sigma_y}{\bar{\sigma}}\right)
\end{equation}
where the parameters are described in Table \ref{tab:plasticity}. A critical aspect of this flow rule is that it is monotonic with respect to $\bar{\sigma}/\sigma_y$, where $\bar{\sigma}=|\boldsymbol{\sigma_0}|$ is the magnitude of the deviatoric shear stress $\boldsymbol{\sigma_0}=\boldsymbol{\sigma}=\frac{1}{3}\mathbf{1}\text{tr}(\boldsymbol{\sigma})$ and $\sigma_y$ is the yield stress, such that plastic deformation does not occur when $\bar{\sigma}/\sigma_y<1$. The second equation is a heat equation that governs the evolution of the effective temperature as
\begin{equation}\label{eq:33}
    c_0 \dot{\chi}=\dfrac{1}{\sigma_y}(\mathbf{D}^{pl}:\boldsymbol{\sigma}_0)(\chi_\infty-\chi) +\nabla\cdot D_\chi\nabla\chi
\end{equation}
where $D_\chi=l^2\sqrt{\mathbf{D}^{pl}:\mathbf{D}^{pl}}$ is a rate-dependent diffusivity. Again, the parameters are described in Table \ref{tab:plasticity}. The parameters are selected to represent a model bulk metallic glass.
\begin{table}[!ht]
 	\centering
 	\caption{Parameters for the STZ plasticity model for a bulk metallic glass material. The Boltzmann constant $k_B=1.38\times 10^{-23}$ J/K is used to express the quantities $\Delta$ and $e_z$ in terms of temperature.}
 	\begin{tabular}{clll}
 		\hline
 		\textbf{Parameter} & \textbf{Unit} &  \textbf{Value} & \textbf{Description} \\
 		\hline \hline
 		$\sigma_y$ &GPa & 0.85 & Yield stress \\
 		$\tau_0$ & sec. &$10^{-13}$ & Molecular vibration timescale \\
 		$\epsilon_0$& - & 0.3 & Typical local strain at STZ transition \\
 		$\Delta/k_B$ &\textrm{K}& $5000$  & Typical activation temperature\\
 		$\tilde{\Omega}/k_B$ & $\text{\AA}^3$ & 349& Typical activation volume \\
 		$T$ &K& 100 & Bath temperature \\
 		$\chi_\infty$ &K& 3000 & Steady-state effective temperature \\
 		$e_z/k_B$ &K& $21000$ & STZ formation energy \\
 		$c_0$ &-& 0.3 & Plastic work fraction \\
 		$l_\chi$ & \text{\AA} & $10.0$ & Diffusion length-scale\\
 		\hline
 	\end{tabular}
    \label{tab:plasticity}
\end{table}

A numerical scheme has been developed by Rycroft et al.\ \cite{Rycroft2015, boffi2020parallel} to solve these equations. The solver, which implements this numerical scheme in an Eulerian finite-difference code under the quasi-static conditions, developed by Rycroft is used here.

\subsubsection{Random Field Initial Conditions}

The STZ theory implies that the effective temperature in the material varies spatially and that these spatial variations play a critical role in the inelastic response of the material. This spatial variation is characterized by an initial $\chi$ field that defines the initial conditions for the simulation. In this work, we will consider two sets of initial conditions for our model material. The first considers $\chi$ as a two-dimensional Gaussian random field having uncertain mean value $\mu_{\chi}$ and coefficient of variation $c_{\chi}$ as described in Table \ref{tab:estochasticity}. Note that the assumption of Gaussianity is based on preliminary results from unpublished prior investigations by Shields and collaborators suggesting that $\chi$ can be reasonable approximated as Gaussian when derived through a coarse-graining procedure from molecular dynamics simulations. Also note that we consider only the parameters of the random field to be uncertain and do not consider uncertainties associated with different realizations of the random $\chi$ field. The purpose of this is to illustrate, through a low-dimensional example, the influence of the relative scale of fluctuations on material behavior. 
\begin{table}[!htb]
	\centering
	\caption{Probability distributions of the STZ random field parameters (2 rvs).}
	\begin{tabular}{cccl}
		\hline
		\textbf{Parameter} & \textbf{Distribution} & \textbf{Range} & \textbf{Description}\\
		\hline \hline
		$\mu_\chi$ & Uniform & [500,700] & Mean of $\chi$ field\\
		$c_\chi$ & Uniform & [0,0.1] & COV of $\chi$ field \\
		\hline
	\end{tabular}
	\label{tab:estochasticity}
\end{table}

\begin{figure}[!htb]
	\centering
	\begin{subfigure}[t]{1.0\textwidth}
	    \includegraphics[width=1.0\textwidth]{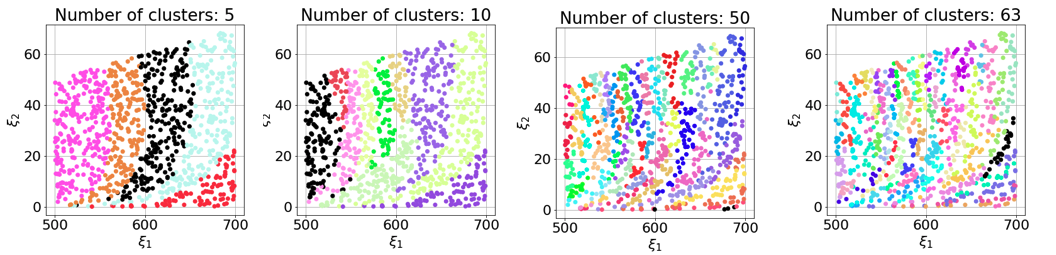}
	    \caption{}
	\end{subfigure}
	\begin{subfigure}[t]{1.0\textwidth}
	    \includegraphics[width=1.0\textwidth]{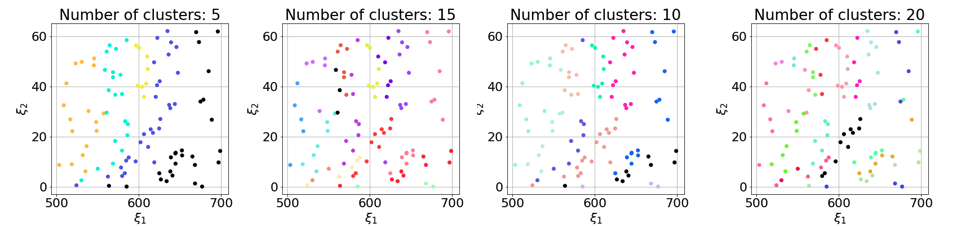}
	    \caption{}
	\end{subfigure}
	\caption{Example 2 -- STZ model with random field initial conditions: Training points clustered according to solution similarity (a) Large data set consisting of 1000 training points (b) Small data set consisting of 100 training points.}
	\label{fig:STZ_clusters}
\end{figure}

To construct the Grassmannian surrogate model, we generate two sets of Monte Carlo training points: a large set containing 1000 samples of the random input vector $\boldsymbol{\xi}_i=(\xi_1, \xi_2)=(\mu_\chi, c_{\chi})$, and a small set containing only 100 samples. The two datasets are shown for increasing number of solution clusters in Figure \ref{fig:STZ_clusters}.
For each realization $i$ of the input parameter vector $\boldsymbol{\xi}_i=(\xi_1, \xi_2)$, the corresponding response is the final plastic strain field after imposing simple shear up to 50\% strain to a simulation box of size $400\text{\AA}\times400\text{\AA}$ having a $32\times 32$ discretization with element size $12.5\text{\AA}\times12.5\text{\AA}$. Hence, the solution matrix $\vec{F}$ is of size $32\times32$. From Figure \ref{fig:STZ_clusters}, we see curved bands of similar solutions for a small number of clusters that give way to a more complex set of bands with increasing number of clusters. We further notice that the solution clusters do not form disjoint regions and therefore do not require parameter space clustering.

\begin{figure}[!htb]
	\centering
	\begin{subfigure}[t]{0.4\textwidth}
		\centering
		\includegraphics[width=\textwidth]{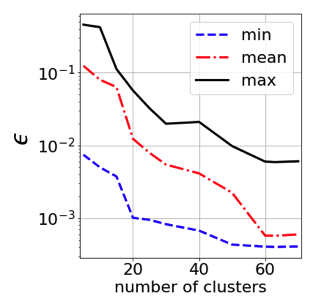}
		\caption{}
	\end{subfigure}%
	~ 
	\begin{subfigure}[t]{0.4\textwidth}
		\centering
		\includegraphics[width=\textwidth]{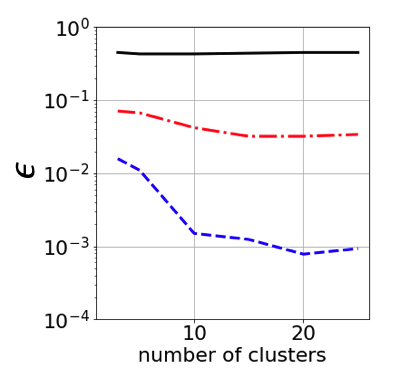}
		\caption{}
	\end{subfigure}
	\caption{Example 2 -- STZ model with random field initial conditions: Minimum, mean and maximum projection error as a function of the number of solution clusters for (a) 1000  and (b) 100 training points.}
	\label{fig:STZ_error1}
\end{figure}

Using the proposed optimization procedure, we find the optimal number of solution clusters to be $n_c=63$ for the large training data set, corresponding to 95\% of the clusters having an error lower than a threshold value of $10^{-3}$ for the error metric in Eq.(\ref{error_mse}). For the small training set, the optimal number of clusters is $n_c=20$ for a relaxed error threshold of $5\times 10^{-2}$. Figure \ref{fig:STZ_error1} shows the minimum, mean, and maximum projection errors used in the optimization from Eq.\ \eqref{error_mse} for increasing number of clusters.  
We notice from these figures that the mean and minimum errors decrease with increasing number of clusters only up to a point. For the small data set, the maximum error does not decrease appreciably with increased number of clusters. This is due to the fact that, with sparse sampling, there are certain regions of the parameter space that simply cannot be improved by modifying the clustering. Only additional training points will improve the solution.

\begin{figure}[!htb]
	\centering
	\includegraphics[width=0.5\textwidth]{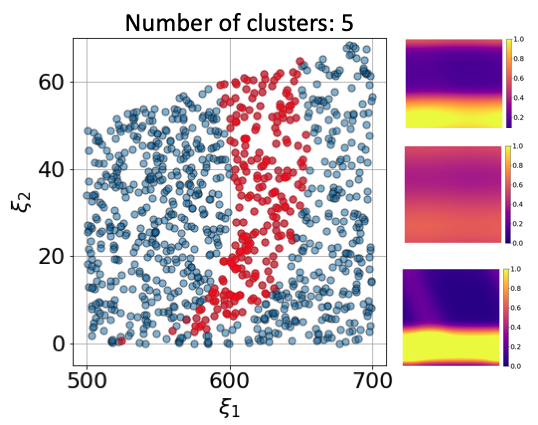}
	\caption{Example 2 -- STZ model with random field initial conditions: Snapshots of the solution at three points inside of a large cluster illustrating that the solution can vary considerably inside clusters that are not selected appropriately.}
	\label{fig:STZ_1000_nc}
\end{figure}

To illustrate the importance of assigning an appropriate number of clusters, consider the cluster shown in Figure \ref{fig:STZ_1000_nc}. This cluster, which is one of only five, has solutions that differ considerably from one another as illustrated by the three snapshots to the right of the training points corresponding to the solution at different points in this cluster.  
Clusters with such varied solutions will lead to errors in the tangent space projection (they are far apart on the Grassmannian) and therefore do not serve as good training sets for the GP.

To train each GP for this model, we again select the Gaussian kernel with initial scale parameter $l=1.0$.  To quantify the performance of the surrogate, we generated $N_{test}=3000$ additional realizations of the random parameter vector $\boldsymbol{\xi}$ and estimated the error metric in Eq.\ \eqref{error_metric} to compare the true solution and the GP surrogate prediction. The GP error is shown as a function of the number of clusters in Figure \ref{fig:STZ_error2}.
\begin{figure}[!htb]
	\centering
	\includegraphics[width=0.40\textwidth]{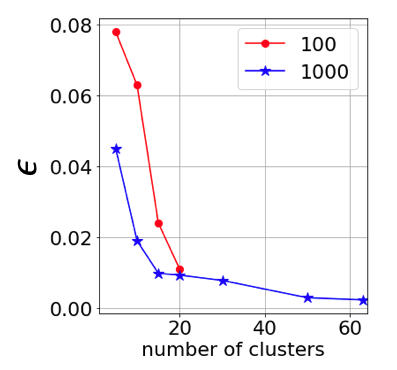}
	\caption{Example 2 -- STZ model with random field initial conditions: GP error  metric of Eq.\ \eqref{error_metric} for increasing number of clusters.}
	\label{fig:STZ_error2}
\end{figure}
To visualize these errors, consider the GP predictions for the ten test points shown in Figure \ref{fig:STZ_1000_testPoints} for different number of solution clusters.
\begin{figure}[!htb]
	\centering
	\includegraphics[width=0.8\textwidth]{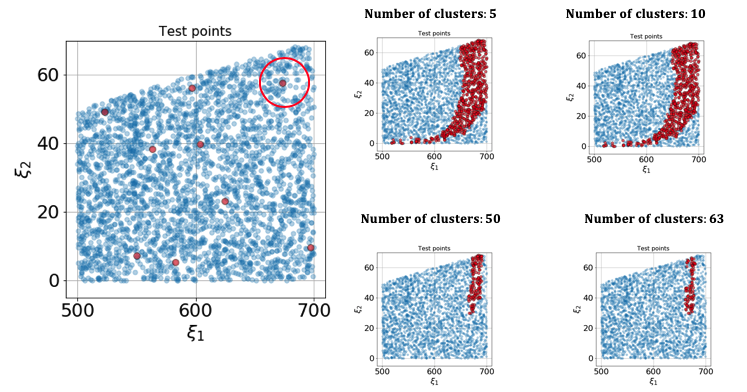}
	\caption{Example 2 -- STZ model with random field initial conditions: Nine test points used to visualize surrogate modeling error. The right images show the clusters to which the circled red test point belongs. }
	\label{fig:STZ_1000_testPoints}
\end{figure}
Figures \ref{fig:STZ_1000_predictions} and \ref{fig:STZ_100_predictions} show the true solution along with the GP predicted solution at each of these 9 test points for different numbers of clusters using the large training set and small training set, respectively.
\begin{figure}[!htb]
	\centering
	\includegraphics[width=1.0\textwidth]{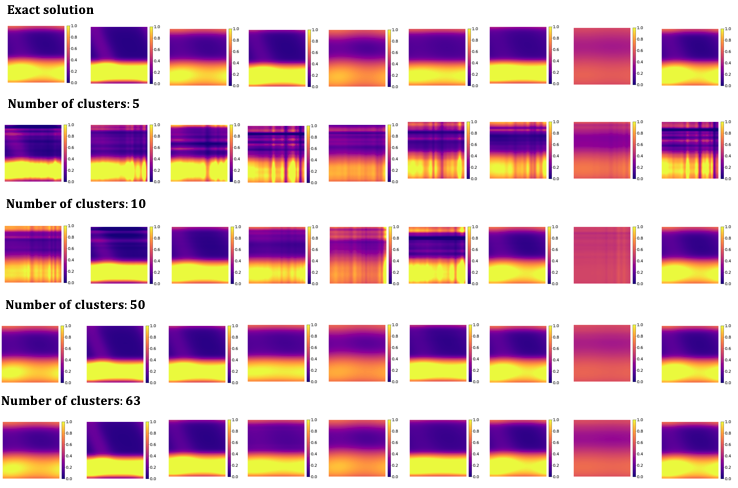}
	\caption{Example 2 -- STZ model with random field initial conditions: Exact solutions (first row) and corresponding predictions for the 9 test points in Figure \ref{fig:STZ_1000_testPoints} for different number of solution clusters with 1000 training points.}
	\label{fig:STZ_1000_predictions}
\end{figure}
\begin{figure}[!htb]
	\centering
	\includegraphics[width=1.0\textwidth]{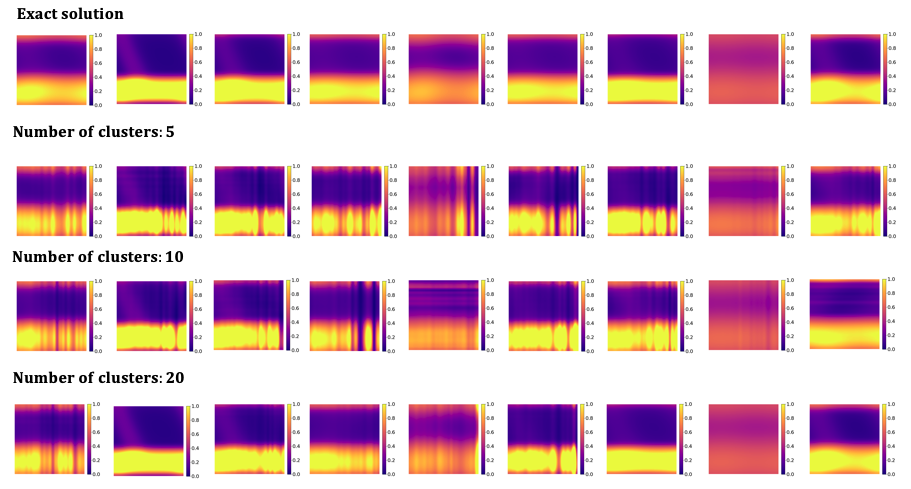}
	\caption{Example 2 -- STZ model with random field initial conditions:  Exact solutions (first row) and corresponding predictions for the 9 test points in Figure \ref{fig:STZ_1000_testPoints} for different number of solution clusters and 100 training points.}
	\label{fig:STZ_100_predictions}
\end{figure}
Here, we again observe that the a small number of clusters results in poor predictions of the solution for many of the points, given the large distances on the manifold that each cluster spans. However, for increasing number of clusters (up to the optimal), the solution predictions improve considerably. The case with only 100 training data has some notable inaccuracies even at the optimal number of clusters, however, due to the fact that the training set simply does not sufficiently span the clusters of material performance. More training data are needed.



\begin{figure}[!htb]
	\centering
	\begin{subfigure}[t]{0.5\textwidth}
		\centering
		\includegraphics[width=\textwidth]{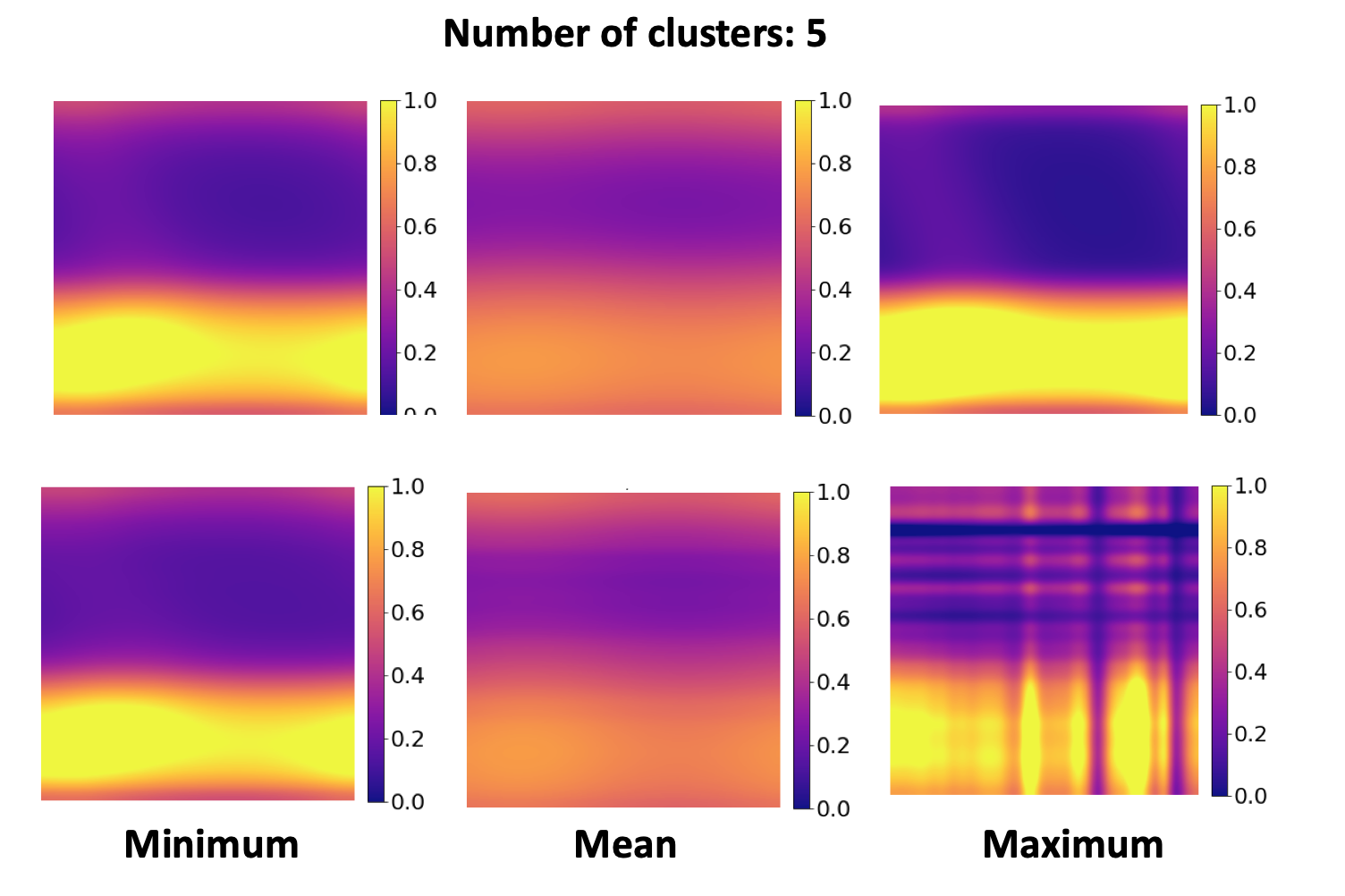}
		\caption{}
	\end{subfigure}%
	~ 
	\begin{subfigure}[t]{0.49\textwidth}
		\centering
		\includegraphics[width=\textwidth]{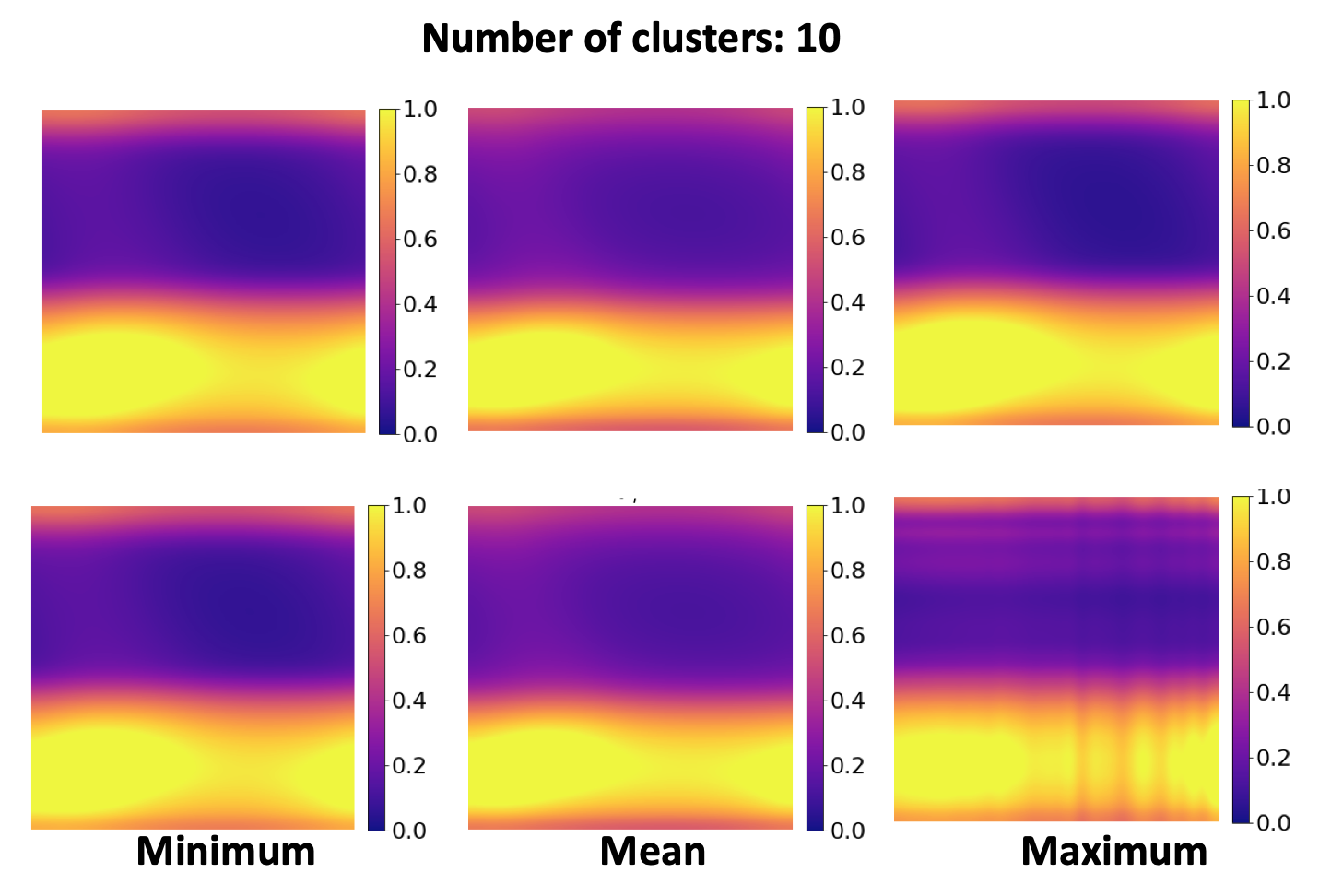}
		\caption{}
	\end{subfigure}

	\begin{subfigure}[t]{0.5\textwidth}
	\centering
	\includegraphics[width=\textwidth]{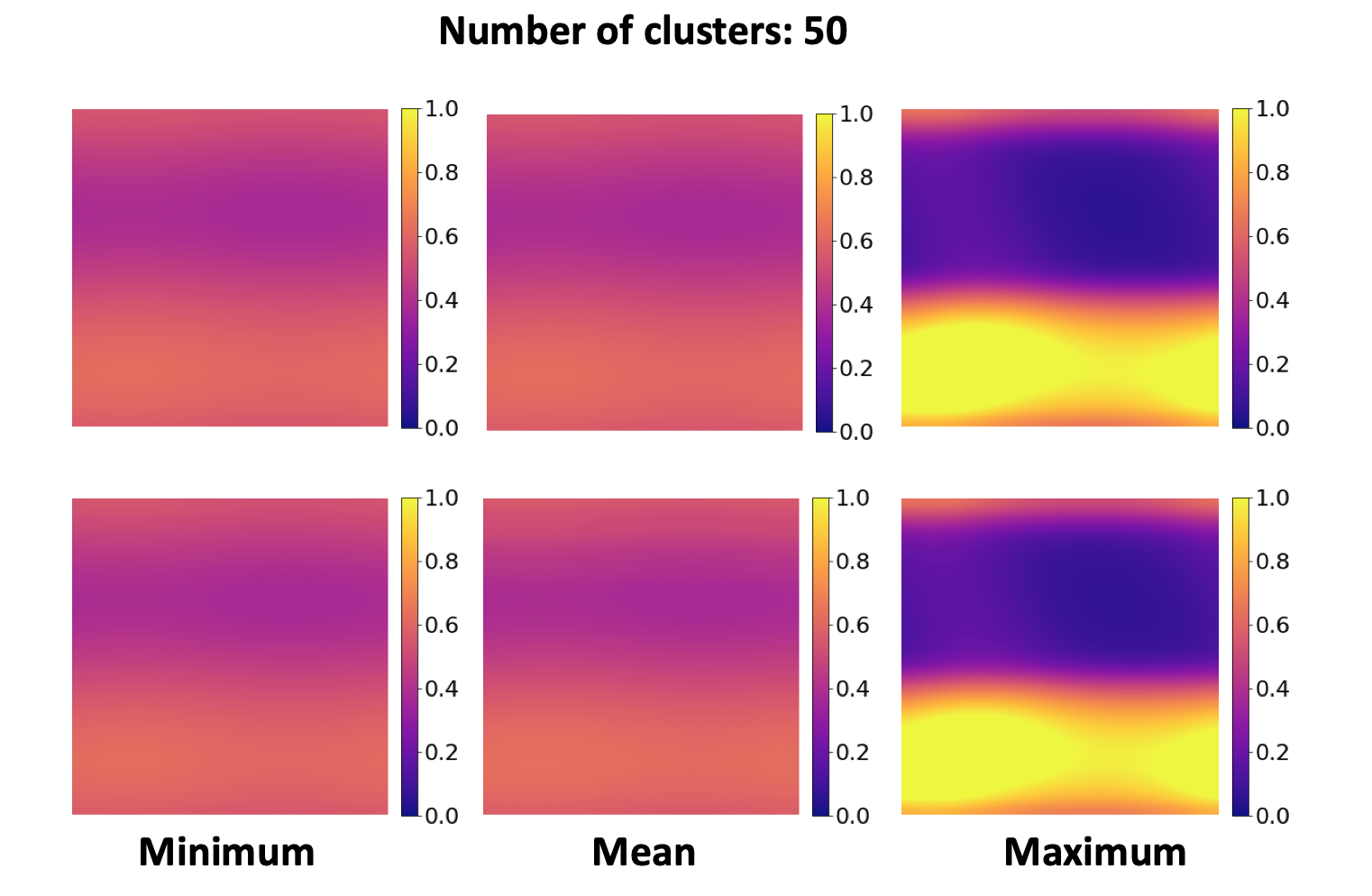}
	\caption{}
\end{subfigure}%
~ 
\begin{subfigure}[t]{0.5\textwidth}
	\centering
	\includegraphics[width=\textwidth]{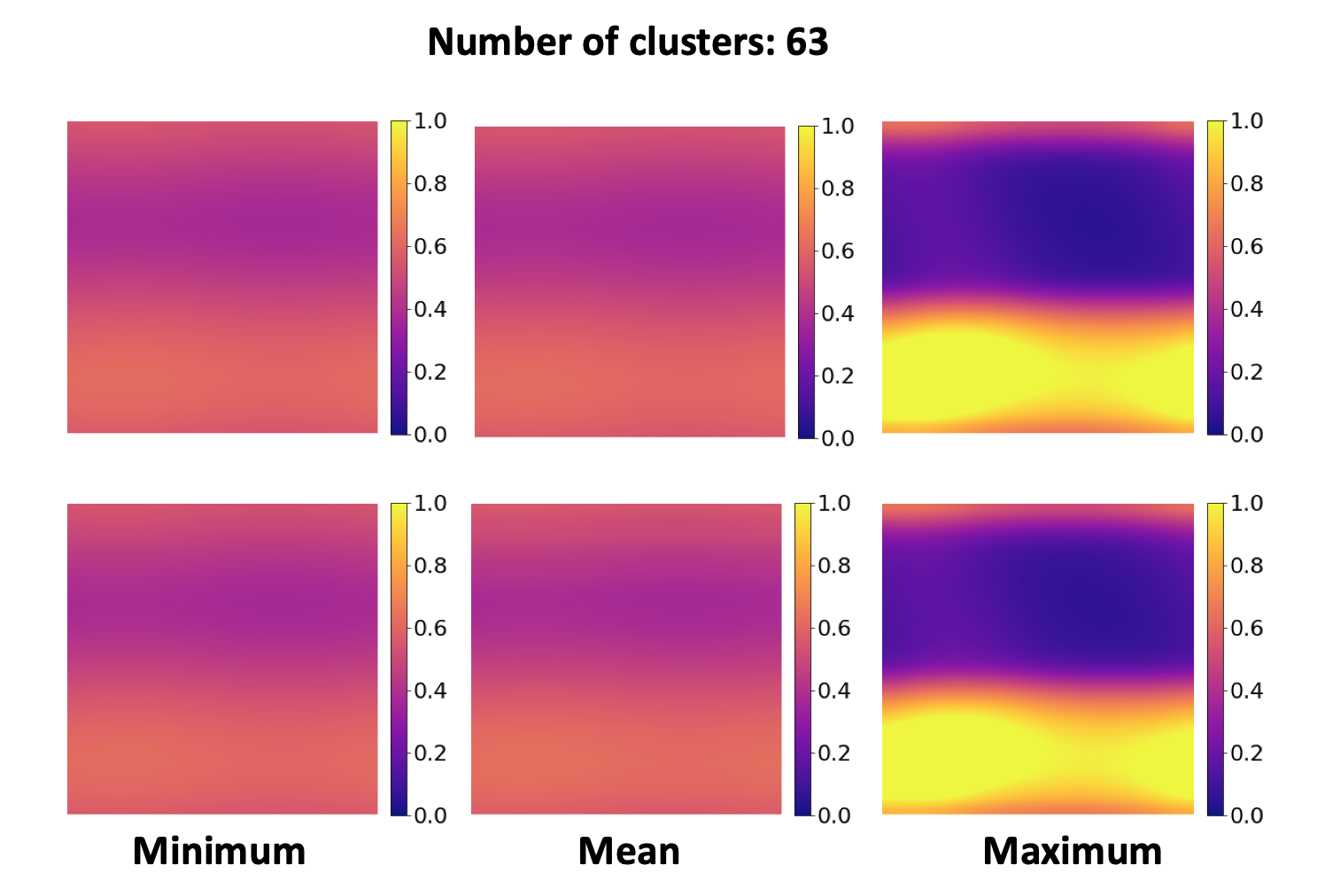}
	\caption{}
\end{subfigure}
	\caption{STZ model -- STZ model with random field initial conditions: Exact solutions (first row) and corresponding  predictions (second row) for the test point that has the minimum  (first column), the mean (second column) and the maximum error (third column)  $\epsilon$ for (a) $n_c=5$, (b) $n_c=10$, (c) $n_c=50$ and (d) $n_c=63$ for 1000 training points.}
	\label{fig:STZ_1000_min_max}
\end{figure}

\begin{figure}[!htb]
	\centering
	\includegraphics[width=0.6\textwidth]{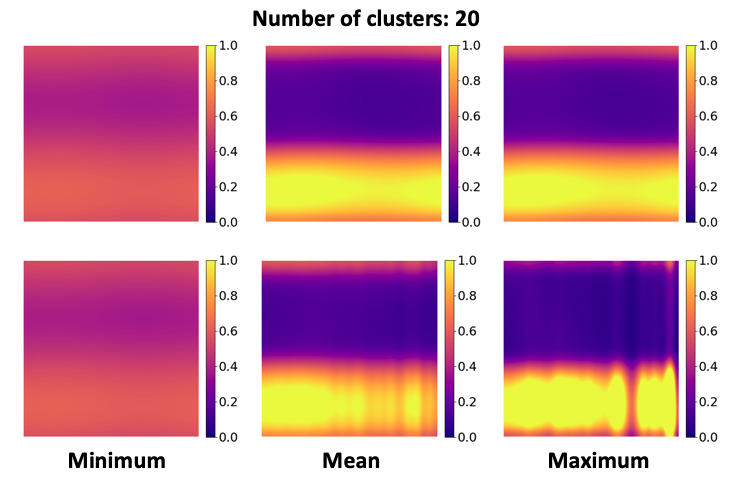}
	\caption{Example 2 -- STZ model with random field initial conditions: Exact solutions (first row) and the corresponding predictions (second row) for the test point that has the minimum  (first column), the mean (second column) and the maximum error (third column)  $\epsilon$ for 100 training points and $n_c=20$.}
	\label{fig:STZ_100_min_max}
\end{figure}

To further illustrate the performance of the  proposed method, Figure \ref{fig:STZ_1000_min_max} shows the exact vs.\ surrogate response from 1000 training points at the test points with the minimum, mean and maximum error for increasing number of solution clusters ($n_c=5, 10, 50, 63$). In all cases, the GP approximation at the points with the minimum and mean error is accurate. However, for the points with maximum error the surrogate approximation is highly dependent on the number of solution clusters. For an insufficient number ($n_c=5, 10$), we can see that the GP fails to capture the true behavior of the model due to the large distance from this solution to the Karcher mean where projection to the tangent space occurs. By increasing the number of solution clusters, we improve the ability of the surrogate to make accurate predictions across the entire parameter space. 
We see the same trend in the case of 100 training points. However, given such a small training set, the surrogate model does not perform well for many of the points even when the optimal number of clusters is selected. This is shown in Figure \ref{fig:STZ_100_min_max}, where we see that the mean and maximum error surrogate solutions have significant error.

\subsection{Initial Conditions Coarse-grained from Molecular Dynamics}

In \cite{Rycroft2017}, a method is proposed for defining initial conditions for STZ simulations by coarse-graining atomistic data from a molecular dynamics simulation of an amorphous material. The coarse graining procedure involves passing a Gaussian window over the atomistic potential energies to produce a weighted average potential energy over a region of atoms centered at grid point $\alpha$ as follows:
\begin{equation}
E_{\alpha} = \frac{\sum_n g_{n} E_n }{\sum_n g_{n}}
\label{eqn:window}
\end{equation}
where $E_n$ is the potential energy of the $n^{th}$ atom and $g_n$ is a Gaussian weight for the $n^{th}$ atom defined as
\begin{equation}
    g_n = \dfrac{2}{\sqrt{2\pi}c}\exp\left(\dfrac{r_n^2}{2c^2}\right)
\end{equation}
where $r_n$ is the distance from atom $n$ to grid point $\alpha$ and $c$ is a length scale that is explored in detail in \cite{Rycroft2017}. A value of $c=50 \text{\AA}$ is used per \cite{Rycroft2017}. The potential energies are then passed through an affine transformation to estimate the effective temperature $\chi_\alpha$ at grid point $\alpha$ as
\begin{equation}
	\chi_{\alpha} = \beta\frac{\sum_n g_n(E_n-E^0)}{\sum_n g_n}
	\label{eqn:affine}
\end{equation}
where $\beta$ and $E^0$ are a material specific constant and a reference energy that corresponds to a state of no disorder, respectively. Notice that Eq.\ \eqref{eqn:affine} incorporates the averaging from Eq.\ \eqref{eqn:window} and the steady-state effective temperature $\chi_\infty$ is calculated as the maximum value $\max(\chi_{\alpha})$. This affine transformation derives from the work of Shi et al.\ \cite{Shi2007}.

In this example, we consider a model 50-50 copper-zirconium metallic glass (CuZr) whose initial $\chi$ field is coarse-grained from an MD simulation. The model, however, has several uncertainties associated with both the coarse-graining parameters ($\beta$ and $E^0$) as well as the parameters of the STZ model as described in Table \ref{tab:plasticity2}.
\begin{table}[!ht]
	\centering
	\caption{Probability distributions of the coarse-graining parameters (top) and STZ plasticity model parameters.}
	\begin{tabular}{lccl}
		\hline
		\textbf{Parameter} & \textbf{Distribution} &  \textbf{Range} & \textbf{Description} \\
		\hline \hline
		$\beta$ &Uniform & [7, 9]& Material specific coarse-graining constant \\
		$E^0$ & Uniform & [-3.37, -3.35] & Reference energy\\ \midrule
		$\Delta/k_B$ & Uniform & [7500, 8500] & Activation temperature\\
		$T$ & Uniform & [50, 150] &Bath temperature\\
		$\tilde{\Omega}/k_B$ & Uniform & [350, 380] & Activation volume \\
		$\epsilon_0$ &Uniform& [0.25, 0.35] & Typical local strain at STZ transition\\
		$c_0$ &Uniform& [0.25, 0.35] & Plastic work fraction \\ 
		$l_{\chi}$& Uniform & [9, 11] & Diffusion length-scale \\
		\hline
	\end{tabular}
	\label{tab:plasticity2}
\end{table}
These parameters are treated here as uncertain because the process of identifying the appropriate parameters for a given glassy material system is an ongoing research field. In \cite{Rycroft2017}, the authors propose an ad hoc means of identifying the appropriate parameters, but more formal procedures remain under investigation. Thus, given the present state of the field, we assign ranges to each parameter and assume a uniform distribution to understand the range of possible material performance given uncertainty in the coarse-graining and the model. This yields a total of eight random variables for our surrogate model.

Again, we apply this for a two-dimensional square grid of $32 \times 32$ equally spaced continuum points (a total of 1024 degrees of freedom) for a simulation of size $400\text{\AA}\times400\text{\AA}$ with 50\% simple shear imposed. The Grassmannian GPs are trained to predict the final plastic strain field at all degrees of freedom for any set of coarse-graining and STZ model parameters in the ranges from Table \ref{tab:plasticity2}. To train the GPs, we execute the STZ model for $N_s=1000$ MCS realizations of the 8-dimensional parameter vector $\{\boldsymbol{\xi}\}$. Applying the optimization procedure, we estimate the optimum number of solutions clusters to be $n_{cl}=61$, corresponding to a threshold value of $10^{-3}$ for the error metric in Eq.\ \eqref{error_mse} for 95\% of the clusters.

Figure \ref{fig:STZ_error2} shows the training error (Eq.\ \eqref{error_mse}) and the GP error (Eq.\ \eqref{error_metric}) for increasing number of clusters. To illustrate the accuracy of the surrogate model predictions, Figure \ref{fig:STZ_8rvs_1200_test} shows the true solution and the Grassmannian GP prediction for 9 test points using the optimal clustering. 
\begin{figure}[!htb]
	\centering
	\begin{subfigure}[t]{0.4\textwidth}
		\centering
		\includegraphics[width=\textwidth]{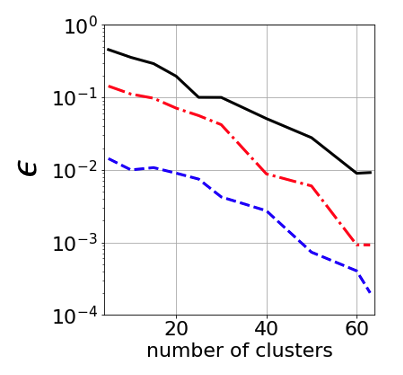}
		\caption{}
	\end{subfigure}%
	~ 
	\begin{subfigure}[t]{0.4\textwidth}
		\centering
		\includegraphics[width=0.98\textwidth]{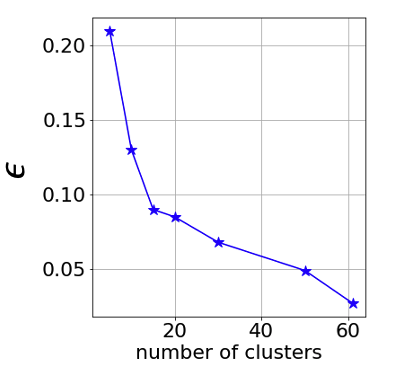}
		\caption{}
	\end{subfigure}
	\caption{Example 2 -- Coarse-grained STZ model: (a) Training error of Eq.\ \eqref{error_mse} and (b) GP error metric of Eq.\ \eqref{error_metric} for increasing number of clusters.}
	\label{fig:STZ_error2}
\end{figure}	
\begin{figure}[!htb]
	\centering
	\includegraphics[width=1.0\textwidth]{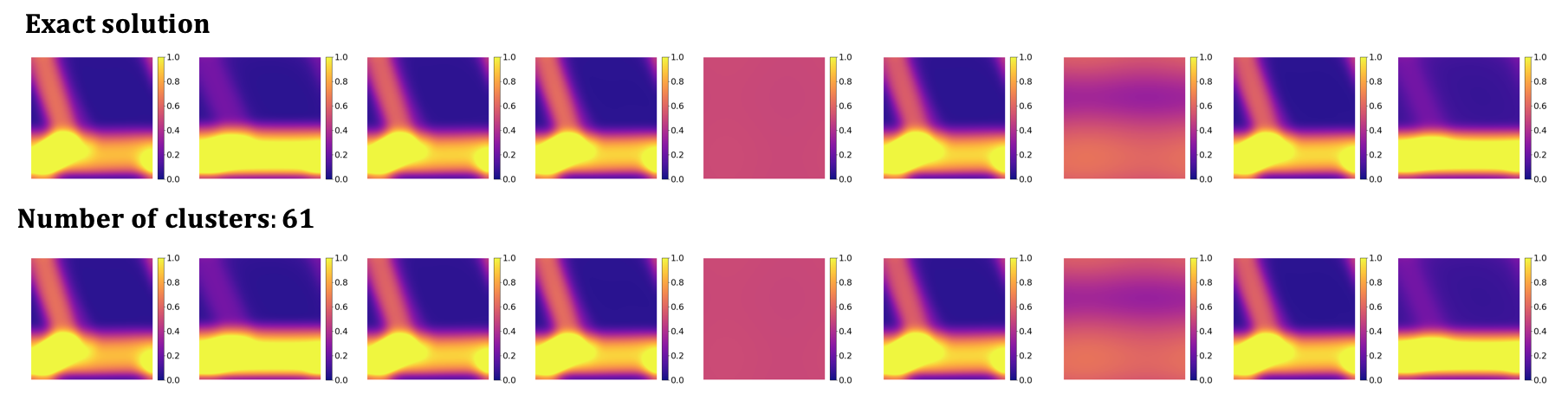}
	\caption{Example 2 -- Coarse-grained STZ model: Exact solutions (first row) and corresponding Grassmannian GP predictions for 9 test points for $n_c=61$ and 1000 training points.}
	\label{fig:STZ_8rvs_1200_test}
\end{figure}
Finally, Figure \ref{fig:STZ_8rvs_min_max} shows  the prediction of the GP at the test points that have the minimum, mean and maximum error, respectively, for the optimum number of solution clusters. From this figure one can see that, at all three points their is a very close agreement between the surrogate prediction and the exact solution. For the minimum ,mean and maximum points the corresponding errors are $3.10\times10^{-4}$, $8.54\times10^{-4}$ and $2.93\times10^{-3}$, respectively.
\begin{figure}[!htb]
	\centering
	\includegraphics[width=0.7\textwidth]{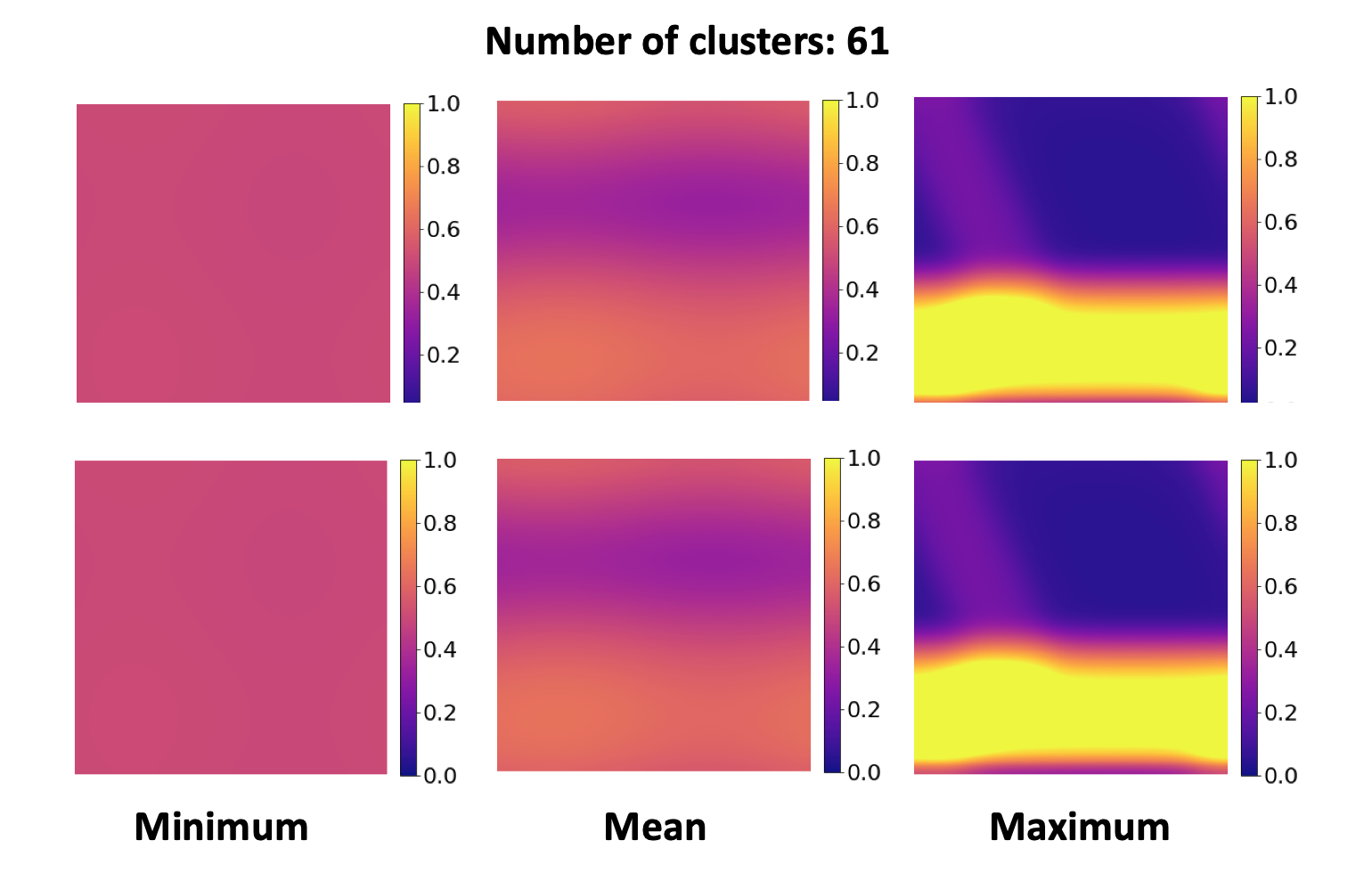}
	\caption{Example 2 -- Coarse-grained STZ model: Exact solution (first row) and the GP prediction (second row) for the test point that has minimum (first column), mean (second column) and maximum error (third column) for 1000 training points and $n_c=61$. }
	\label{fig:STZ_8rvs_min_max}
\end{figure}

\section{Additional computational cost benefits}

\noindent

One computational bottleneck of the standard GP is to compute the Cholesky decomposition of the covariance function $k(\cdot, \cdot)$ required for the estimation of the hyper-parameters of the GP model for large data sets \cite{wang2020accelerated}. For a training set containing $N_s$ samples and a scalar output the complexity of training is $\mathcal{O}(N_s^3)$ and that of making a prediction is of $\mathcal{O}(N_s)$ for the mean and $\mathcal{O}(N_s^2)$ for the variance. Moreover, when we are dealing with multiple outputs (e.g. vectors, matrices) the training time of the GP increases proportional to the size $n=m_f\times n_f$ of the output and becomes  $n\times \mathcal{O}(N_s^3)$.


Comparing the computational training time of the standard GP with the proposed local GPs, we achieve a significant reduction in the computational cost given (a) the much smaller number of training points for each local surrogate and (b) the reduced output space. For  each  solution cluster $i$ for $i=1, \ldots, n_{cl}$  the local GP is linked to a covariance function of dimension $(N_i\times N_i)$ where  $N_i$ is the number of training samples in cluster $i$. Moreover, each output has a reduced dimension of $m_f\times r$. Since $N_i <<N_s$ and $r<<n_f$ the computational cost of the GP training is orders of magnitude lower than the standard GP. A direct comparison of the training time between the proposed local GPs and a full GP trained on the plastic strain fields themselves for the STZ model with 2 random variables with 1000 training points, illustrates that the proposed method is far more efficient and yields a more accurate solution. The full GP requires $t=2,897$ sec.\ for training with an average error $\epsilon=0.023$. Meanwhile, the total training time for the local Grassmannian  GPs for $n_{cl}= 63$ clusters is $t=293$ sec.\ with corresponding error $\epsilon=0.0032$.


\section*{Conclusions}

This paper introduces new machine learning methods for surrogate model construction in the context of data-driven uncertainty quantification (UQ) for high-dimensional (input \& output) complex physical/engineering systems. We illustrate how the solution of a high-dimensional model can be represented on a lower-dimensional Riemannian manifold (specifically the Grassmannian), which enables the construction of surrogate models for solution prediction. In particular, the manifold representation allows areas of the parameter space to be identified where the model has similar behavior. We then propose a ``solution clustering'' scheme and build local Gaussian process regression models that map points in the parameter space to the reduced solution in the tangent space of the Grassmannian to predict the solution at new points in the parameter space.

The proposed method is outlined in Figure {\color{red}} and involves five steps: 1. Evaluation the model of a set of training points and projecting them onto the Grassmann manifold; 2. Solution clustering using a manifold clustering technique; 3. A mapping of each cluster to the tangents space of the Grassmannian at the Karcher mean of each cluster; 4. An optional parameter space clustering that is necessary only if solution clusters form disjoint clusters in the parameter space; and 5. Construction of local Gaussian Process models for each cluster for solution preduction.


We demonstrate the proposed method using two examples. The first is the Kraichnan-Orszag (KO) three mode problem, which is a nonlinear three-dimensional set of stochastic ordinary differential equations subjected to random initial conditions. We show that a two-level clustering (solution \& input) is required in order to accurately train local GP surrogates and predict the time history response for realizations of the initial conditions. The second problem considers the plastic deformation of an amorphous solid using the shear transformation zone (STZ) theory of plasticity. We consider two cases for this model. The first is a low-dimensional example that is used to illustrate the solution clustering scheme and considers only the mean and coefficient of variation of an initial random $\chi$ field as the uncertain parameters. The second case is a more practically relevant case that considers uncertainty in the parameters used to coarse-grain a molecular dynamics simulation for STZ modeling along with uncertainty in the STZ model parameters themselves. In total, this case has eight random parameters. In both cases, we train local surrogates having similar solutions and that provide very accurate predictions of the material response.


\section*{Acknowledgements}

\noindent
This material is based upon work supported by the U.S. Department of Energy, Office of Science, Office of Advanced Scientific Computing Research Early Career Research program under Award Number DE-SC0020428 with Dr. Steven Lee as program manager.





 \bibliographystyle{model1-num-names}
 \bibliography{References.bib}

\begin{thebibliography}{72}
\expandafter\ifx\csname natexlab\endcsname\relax\def\natexlab#1{#1}\fi
\providecommand{\bibinfo}[2]{#2}
\ifx\xfnm\relax \def\xfnm[#1]{\unskip,\space#1}\fi
\bibitem[{Xiu and Karniadakis(2002)}]{xiu2002wiener}
\bibinfo{author}{D.~Xiu}, \bibinfo{author}{G.~E. Karniadakis},
\newblock \bibinfo{title}{The wiener--askey polynomial chaos for stochastic
  differential equations},
\newblock \bibinfo{journal}{SIAM journal on scientific computing}
  \bibinfo{volume}{24} (\bibinfo{year}{2002}) \bibinfo{pages}{619--644}.
\bibitem[{Wan and Karniadakis(2005)}]{MEgPC2005}
\bibinfo{author}{X.~Wan}, \bibinfo{author}{G.~E. Karniadakis},
\newblock \bibinfo{title}{An adaptive multi-element generalized polynomial
  chaos method for stochastic differential equations},
\newblock \bibinfo{journal}{Journal of Computational Physics}
  \bibinfo{volume}{209} (\bibinfo{year}{2005}) \bibinfo{pages}{617--642}.
\bibitem[{Wan and Karniadakis(2006)}]{MEgCM2006}
\bibinfo{author}{X.~Wan}, \bibinfo{author}{G.~E. Karniadakis},
\newblock \bibinfo{title}{Multi-element generalized polynomial chaos for
  arbitrary probability measures},
\newblock \bibinfo{journal}{SIAM Journal on Scientific Computing}
  \bibinfo{volume}{28} (\bibinfo{year}{2006}) \bibinfo{pages}{901--928}.
\bibitem[{Blatman and Sudret(2011)}]{Blatman2011}
\bibinfo{author}{G.~Blatman}, \bibinfo{author}{B.~Sudret},
\newblock \bibinfo{title}{Adaptive sparse polynomial chaos expansion based on
  least angle regression},
\newblock \bibinfo{journal}{J. Comput. Phys.} \bibinfo{volume}{230}
  (\bibinfo{year}{2011}) \bibinfo{pages}{2345--2367}.
\bibitem[{Papaioannou et~al.(2019)Papaioannou, Ehre, and
  Straub}]{papaioannou2019pls}
\bibinfo{author}{I.~Papaioannou}, \bibinfo{author}{M.~Ehre},
  \bibinfo{author}{D.~Straub},
\newblock \bibinfo{title}{Pls-based adaptation for efficient pce representation
  in high dimensions},
\newblock \bibinfo{journal}{Journal of Computational Physics}
  \bibinfo{volume}{387} (\bibinfo{year}{2019}) \bibinfo{pages}{186--204}.
\bibitem[{L{\"u}then et~al.(2020)L{\"u}then, Marelli, and
  Sudret}]{luthen2020sparse}
\bibinfo{author}{N.~L{\"u}then}, \bibinfo{author}{S.~Marelli},
  \bibinfo{author}{B.~Sudret},
\newblock \bibinfo{title}{Sparse polynomial chaos expansions: Literature survey
  and benchmark},
\newblock \bibinfo{journal}{arXiv preprint arXiv:2002.01290}
  (\bibinfo{year}{2020}).
\bibitem[{Babuška et~al.(2007)Babuška, Nobile, and Tempone}]{Babuska2007}
\bibinfo{author}{I.~Babuška}, \bibinfo{author}{F.~Nobile},
  \bibinfo{author}{R.~Tempone},
\newblock \bibinfo{title}{A stochastic collocation method for elliptic partial
  differential equations with random input data},
\newblock \bibinfo{journal}{SIAM J. Numer. Anal.} \bibinfo{volume}{45}
  (\bibinfo{year}{2007}) \bibinfo{pages}{1005--1034}.
\bibitem[{Foo et~al.(2008)Foo, Wan, and Karniadakis}]{MEPCM2008}
\bibinfo{author}{J.~Foo}, \bibinfo{author}{X.~Wan},
  \bibinfo{author}{G.~Karniadakis},
\newblock \bibinfo{title}{The multi-element probabilistic collocation method
  {(ME-PCM)}: Error analysis and applications},
\newblock \bibinfo{journal}{Journal of Computational Physics}
  \bibinfo{volume}{227} (\bibinfo{year}{2008}) \bibinfo{pages}{9572--9595}.
\bibitem[{Ma and Zabaras(2009)}]{ASGC2009}
\bibinfo{author}{X.~Ma}, \bibinfo{author}{N.~Zabaras},
\newblock \bibinfo{title}{An adaptive hierarchical sparse grid collocation
  algorithm for the solution of stochastic differential equations},
\newblock \bibinfo{journal}{Journal of Computational Physics}
  \bibinfo{volume}{228} (\bibinfo{year}{2009}) \bibinfo{pages}{3084--3113}.
\bibitem[{Xiu and Hesthaven(2005)}]{Xiu2005}
\bibinfo{author}{D.~Xiu}, \bibinfo{author}{J.~S. Hesthaven},
\newblock \bibinfo{title}{High-order collocation methods for differential
  equations with random inputs},
\newblock \bibinfo{journal}{SIAM J. Sci. Comput.} \bibinfo{volume}{27}
  (\bibinfo{year}{2005}) \bibinfo{pages}{1118--1139}.
\bibitem[{Loukrezis et~al.(2019)Loukrezis, R{\"o}mer, and
  De~Gersem}]{loukrezis2019assessing}
\bibinfo{author}{D.~Loukrezis}, \bibinfo{author}{U.~R{\"o}mer},
  \bibinfo{author}{H.~De~Gersem},
\newblock \bibinfo{title}{Assessing the performance of leja and clenshaw-curtis
  collocation for computational electromagnetics with random input data},
\newblock \bibinfo{journal}{International Journal for Uncertainty
  Quantification} \bibinfo{volume}{9} (\bibinfo{year}{2019}).
\bibitem[{Krige(1951)}]{Krige51}
\bibinfo{author}{D.~G. Krige},
\newblock \bibinfo{title}{A statistical approach to some basic mine valuation
  problems on the witwatersrand},
\newblock \bibinfo{journal}{J. Chem. Metall. Mining Soc. South Africa}
  \bibinfo{volume}{52} (\bibinfo{year}{1951}) \bibinfo{pages}{119--139}.
\bibitem[{Santner et~al.(2003)Santner, Williams, and Notz}]{Santner2003}
\bibinfo{author}{T.~J. Santner}, \bibinfo{author}{B.~J. Williams},
  \bibinfo{author}{W.~I. Notz}, \bibinfo{title}{The Design and Analysis of
  Computer Experiments}, \bibinfo{publisher}{Springer}, \bibinfo{year}{2003}.
\bibitem[{Rasmussen and Williams(2006)}]{Rasmussen2006}
\bibinfo{author}{C.~E. Rasmussen}, \bibinfo{author}{C.~K.~I. Williams},
\newblock \bibinfo{title}{Gaussian processes for machine learning. {A}daptive
  computation and machine learning}  (\bibinfo{year}{2006}).
\bibitem[{Echard et~al.(2011)Echard, Gayton, and Lemaire}]{Echard2011}
\bibinfo{author}{B.~Echard}, \bibinfo{author}{N.~Gayton},
  \bibinfo{author}{M.~Lemaire},
\newblock \bibinfo{title}{{AK-MCS}: an active learning reliability method
  combining kriging and monte carlo simulation},
\newblock \bibinfo{journal}{Struct. Saf.} \bibinfo{volume}{33}
  (\bibinfo{year}{2011}) \bibinfo{pages}{145--154}.
\bibitem[{Schobi et~al.(2015)Schobi, Sudret, and Wiart}]{schobi2015polynomial}
\bibinfo{author}{R.~Schobi}, \bibinfo{author}{B.~Sudret},
  \bibinfo{author}{J.~Wiart},
\newblock \bibinfo{title}{Polynomial-chaos-based kriging},
\newblock \bibinfo{journal}{International Journal for Uncertainty
  Quantification} \bibinfo{volume}{5} (\bibinfo{year}{2015}).
\bibitem[{Csat\'{o} and Opper(2002)}]{Csato2002}
\bibinfo{author}{L.~Csat\'{o}}, \bibinfo{author}{M.~Opper},
\newblock \bibinfo{title}{Sparse online gaussian processes},
\newblock \bibinfo{journal}{N.Comp.} \bibinfo{volume}{14}
  (\bibinfo{year}{2002}) \bibinfo{pages}{641--668}.
\bibitem[{Smola and Bartlett(2001)}]{Smola2001}
\bibinfo{author}{A.~Smola}, \bibinfo{author}{P.~Bartlett},
\newblock \bibinfo{title}{Sparse greedy gaussian process regression},
\newblock \bibinfo{journal}{In Advances in NIPS 13}  (\bibinfo{year}{2001})
  \bibinfo{pages}{619--625}.
\bibitem[{Wang et~al.(2020)Wang, Leiter, Plech{\'a}{\v{c}}, and
  Knap}]{wang2020accelerated}
\bibinfo{author}{T.~Wang}, \bibinfo{author}{K.~W. Leiter},
  \bibinfo{author}{P.~Plech{\'a}{\v{c}}}, \bibinfo{author}{J.~Knap},
\newblock \bibinfo{title}{Accelerated scale bridging with sparsely approximated
  gaussian learning},
\newblock \bibinfo{journal}{Journal of Computational Physics}
  \bibinfo{volume}{403} (\bibinfo{year}{2020}) \bibinfo{pages}{109049}.
\bibitem[{Rozza et~al.(2008)Rozza, Huynh, and Patera}]{Rozza2008}
\bibinfo{author}{G.~Rozza}, \bibinfo{author}{D.~Huynh},
  \bibinfo{author}{A.~Patera},
\newblock \bibinfo{title}{Reduced basis approximation and a posteriori error
  estimation for affinely parametrized elliptic coercive partial differential
  equations},
\newblock \bibinfo{journal}{Arch. Comput. Methods Eng.} \bibinfo{volume}{15}
  (\bibinfo{year}{2008}) \bibinfo{pages}{229--275}.
\bibitem[{Tan et~al.(2008)Tan, Willcox, and Ghattas}]{Tan2008}
\bibinfo{author}{B.~Tan}, \bibinfo{author}{K.~Willcox},
  \bibinfo{author}{O.~Ghattas},
\newblock \bibinfo{title}{Model reduction for large-scale systems with
  high-dimensional parametric input space},
\newblock \bibinfo{journal}{SIAM J. Sci. Comput.} \bibinfo{volume}{30}
  (\bibinfo{year}{2008}) \bibinfo{pages}{3270--3288}.
\bibitem[{Benner et~al.(2015)Benner, Gugercin, and Willcox}]{Benner2015}
\bibinfo{author}{P.~Benner}, \bibinfo{author}{S.~Gugercin},
  \bibinfo{author}{K.~Willcox},
\newblock \bibinfo{title}{A survey of projection-based model reduction methods
  for parametric dynamical systems},
\newblock \bibinfo{journal}{SIAM Review} \bibinfo{volume}{57}
  (\bibinfo{year}{2015}) \bibinfo{pages}{483--531}.
\bibitem[{Chatterjee(2000)}]{POD2000}
\bibinfo{author}{A.~Chatterjee},
\newblock \bibinfo{title}{An introduction to the proper orthogonal
  decomposition},
\newblock \bibinfo{journal}{Current Science} \bibinfo{volume}{78}
  (\bibinfo{year}{2000}) \bibinfo{pages}{808--817}.
\bibitem[{Zou et~al.(2019)Zou, Kouri, and Aquino}]{zou2019adaptive}
\bibinfo{author}{Z.~Zou}, \bibinfo{author}{D.~Kouri},
  \bibinfo{author}{W.~Aquino},
\newblock \bibinfo{title}{An adaptive local reduced basis method for solving
  pdes with uncertain inputs and evaluating risk},
\newblock \bibinfo{journal}{Computer Methods in Applied Mechanics and
  Engineering} \bibinfo{volume}{345} (\bibinfo{year}{2019})
  \bibinfo{pages}{302--322}.
\bibitem[{Stabile and Rosic(2020)}]{stabile2020bayesian}
\bibinfo{author}{G.~Stabile}, \bibinfo{author}{B.~Rosic},
\newblock \bibinfo{title}{Bayesian identification of a projection-based reduced
  order model for computational fluid dynamics},
\newblock \bibinfo{journal}{Computers \& Fluids} \bibinfo{volume}{201}
  (\bibinfo{year}{2020}) \bibinfo{pages}{104477}.
\bibitem[{Amsallem and Farhat(2008)}]{Amsallem2008}
\bibinfo{author}{D.~Amsallem}, \bibinfo{author}{C.~Farhat},
\newblock \bibinfo{title}{Interpolation method for adapting reduced-order
  models and application to aeroelasticity},
\newblock \bibinfo{journal}{AIAA J.} \bibinfo{volume}{46}
  (\bibinfo{year}{2008}) \bibinfo{pages}{1803--1813}.
\bibitem[{Amsallem and Farhat(2011)}]{Amsallem2011}
\bibinfo{author}{D.~Amsallem}, \bibinfo{author}{C.~Farhat},
\newblock \bibinfo{title}{An online method for interpolating linear parametric
  reduced-order models},
\newblock \bibinfo{journal}{SIAM J. Sci. Comput.} \bibinfo{volume}{33}
  (\bibinfo{year}{2011}) \bibinfo{pages}{2169--2198}.
\bibitem[{Jolliffe(1986)}]{PCA1986}
\bibinfo{author}{I.~Jolliffe}, \bibinfo{title}{Principal component analysis and
  factor analysis. In: Principal Component Analysis},
  \bibinfo{publisher}{Springer}, \bibinfo{address}{New York},
  \bibinfo{year}{1986}.
\bibitem[{Roweis and Saul(2000)}]{Roweis2000}
\bibinfo{author}{S.~Roweis}, \bibinfo{author}{L.~Saul},
\newblock \bibinfo{title}{Nonlinear dimensionality reduction by locally linear
  embedding},
\newblock \bibinfo{journal}{Science} \bibinfo{volume}{290}
  (\bibinfo{year}{2000}) \bibinfo{pages}{2323--2326}.
\bibitem[{Belkin and Niyogi(2003)}]{Belkin2003}
\bibinfo{author}{M.~Belkin}, \bibinfo{author}{P.~Niyogi},
\newblock \bibinfo{title}{Laplacian eigenmaps for dimensionality reduction and
  data representation},
\newblock \bibinfo{journal}{Neural Computation} \bibinfo{volume}{15}
  (\bibinfo{year}{2003}) \bibinfo{pages}{1373--1396}.
\bibitem[{Coifman and Lafon(2006)}]{DMaps2005}
\bibinfo{author}{R.~Coifman}, \bibinfo{author}{S.~Lafon},
\newblock \bibinfo{title}{Diffusion maps},
\newblock \bibinfo{journal}{Applied and Computational Harmonic Analysis}
  \bibinfo{volume}{21} (\bibinfo{year}{2006}) \bibinfo{pages}{5--30}.
\bibitem[{Tenenbaum et~al.(2000)Tenenbaum, de~Silva, and
  Langford}]{Tenenbaum2000}
\bibinfo{author}{J.~Tenenbaum}, \bibinfo{author}{V.~de~Silva},
  \bibinfo{author}{J.~Langford},
\newblock \bibinfo{title}{A global geometric framework for nonlinear
  dimensionality reduction},
\newblock \bibinfo{journal}{Science} \bibinfo{volume}{290}
  (\bibinfo{year}{2000}) \bibinfo{pages}{2319--2323}.
\bibitem[{Donoho and Grimes(2003)}]{Dohoho2003}
\bibinfo{author}{D.~Donoho}, \bibinfo{author}{C.~Grimes},
\newblock \bibinfo{title}{Hessian eigenmaps: locally linear embedding
  techniques for high-dimensional data},
\newblock \bibinfo{journal}{Proc. Natl. Acad. Sci.} \bibinfo{volume}{100}
  (\bibinfo{year}{2003}) \bibinfo{pages}{5591--5596}.
\bibitem[{Sch\"{o}lkopf et~al.(1998)Sch\"{o}lkopf, Smola, and
  Müller}]{Scholkopf1998}
\bibinfo{author}{B.~Sch\"{o}lkopf}, \bibinfo{author}{A.~Smola},
  \bibinfo{author}{K.~R. Müller},
\newblock \bibinfo{title}{Nonlinear component analysis as a kernel eigenvalue
  problem},
\newblock \bibinfo{journal}{Neural Computation} \bibinfo{volume}{10}
  (\bibinfo{year}{1998}) \bibinfo{pages}{1299--1319}.
\bibitem[{Zhang and Zha(2004)}]{Zhang2004}
\bibinfo{author}{Z.~Zhang}, \bibinfo{author}{H.~Zha},
\newblock \bibinfo{title}{Principal manifolds and nonlinear dimensionality
  reduction via tangent space alignment},
\newblock \bibinfo{journal}{SIAM J Sci Comput.} \bibinfo{volume}{26}
  (\bibinfo{year}{2004}) \bibinfo{pages}{313--338}.
\bibitem[{Lataniotis et~al.(2020)Lataniotis, Marelli, and
  Sudret}]{lataniotis2020extending}
\bibinfo{author}{C.~Lataniotis}, \bibinfo{author}{S.~Marelli},
  \bibinfo{author}{B.~Sudret},
\newblock \bibinfo{title}{Extending classical surrogate modeling to high
  dimensions through supervised dimensionality reduction: A data-driven
  approach},
\newblock \bibinfo{journal}{International Journal for Uncertainty
  Quantification} \bibinfo{volume}{10} (\bibinfo{year}{2020}).
\bibitem[{Raissi et~al.(2019)Raissi, Perdikaris, and
  Karniadakis}]{raissi2019physics}
\bibinfo{author}{M.~Raissi}, \bibinfo{author}{P.~Perdikaris},
  \bibinfo{author}{G.~E. Karniadakis},
\newblock \bibinfo{title}{Physics-informed neural networks: A deep learning
  framework for solving forward and inverse problems involving nonlinear
  partial differential equations},
\newblock \bibinfo{journal}{Journal of Computational Physics}
  \bibinfo{volume}{378} (\bibinfo{year}{2019}) \bibinfo{pages}{686--707}.
\bibitem[{Wang et~al.(2020)Wang, Broccardo, and Song}]{wang2020probabilistic}
\bibinfo{author}{Z.~Wang}, \bibinfo{author}{M.~Broccardo},
  \bibinfo{author}{J.~Song},
\newblock \bibinfo{title}{Probabilistic performance-pattern decomposition
  (pppd): analysis framework and applications to stochastic mechanical
  systems},
\newblock \bibinfo{journal}{arXiv preprint arXiv:2003.02205}
  (\bibinfo{year}{2020}).
\bibitem[{Kalogeris and Papadopoulos(2019)}]{kalogeris2019diffusion}
\bibinfo{author}{I.~Kalogeris}, \bibinfo{author}{V.~Papadopoulos},
\newblock \bibinfo{title}{Diffusion maps-based surrogate modeling: An
  alternative machine learning approach},
\newblock \bibinfo{journal}{International Journal for Numerical Methods in
  Engineering}  (\bibinfo{year}{2019}).
\bibitem[{Soize and Ghanem(2016)}]{soize2016data}
\bibinfo{author}{C.~Soize}, \bibinfo{author}{R.~Ghanem},
\newblock \bibinfo{title}{Data-driven probability concentration and sampling on
  manifold},
\newblock \bibinfo{journal}{Journal of Computational Physics}
  \bibinfo{volume}{321} (\bibinfo{year}{2016}) \bibinfo{pages}{242--258}.
\bibitem[{Soize and Farhat(2017)}]{soize2017nonparametric}
\bibinfo{author}{C.~Soize}, \bibinfo{author}{C.~Farhat},
\newblock \bibinfo{title}{A nonparametric probabilistic approach for
  quantifying uncertainties in low-dimensional and high-dimensional nonlinear
  models},
\newblock \bibinfo{journal}{International Journal for Numerical Methods in
  Engineering} \bibinfo{volume}{109} (\bibinfo{year}{2017})
  \bibinfo{pages}{837--888}.
\bibitem[{Farhat et~al.(2018)Farhat, Bos, Avery, and
  Soize}]{farhat2018modeling}
\bibinfo{author}{C.~Farhat}, \bibinfo{author}{A.~Bos},
  \bibinfo{author}{P.~Avery}, \bibinfo{author}{C.~Soize},
\newblock \bibinfo{title}{Modeling and quantification of model-form
  uncertainties in eigenvalue computations using a stochastic reduced model},
\newblock \bibinfo{journal}{AIAA Journal} \bibinfo{volume}{56}
  (\bibinfo{year}{2018}) \bibinfo{pages}{1198--1210}.
\bibitem[{Soize and Farhat(2019)}]{soize2019probabilistic}
\bibinfo{author}{C.~Soize}, \bibinfo{author}{C.~Farhat},
\newblock \bibinfo{title}{Probabilistic learning for modeling and quantifying
  model-form uncertainties in nonlinear computational mechanics},
\newblock \bibinfo{journal}{International Journal for Numerical Methods in
  Engineering} \bibinfo{volume}{117} (\bibinfo{year}{2019})
  \bibinfo{pages}{819--843}.
\bibitem[{Giovanis and Shields(2018)}]{Giovanis2018}
\bibinfo{author}{D.~Giovanis}, \bibinfo{author}{M.~Shields},
\newblock \bibinfo{title}{Uncertainty quantification for complex systems with
  very high dimensional response using grassmann manifold variations},
\newblock \bibinfo{journal}{J Comput Phys.} \bibinfo{volume}{364}
  (\bibinfo{year}{2018}) \bibinfo{pages}{393--415}.
\bibitem[{Giovanis and Shields(2019)}]{Giovanis2019}
\bibinfo{author}{D.~Giovanis}, \bibinfo{author}{M.~Shields},
\newblock \bibinfo{title}{Variance-based simplex stochastic collocation with
  model order reduction for high-dimensional systems},
\newblock \bibinfo{journal}{J. Numer. Methods Eng.} \bibinfo{volume}{117}
  (\bibinfo{year}{2019}) \bibinfo{pages}{1079--1116}.
\bibitem[{Wang et~al.(2019)Wang, Guilleminot, and Soize}]{wang2019modeling}
\bibinfo{author}{H.~Wang}, \bibinfo{author}{J.~Guilleminot},
  \bibinfo{author}{C.~Soize},
\newblock \bibinfo{title}{Modeling uncertainties in molecular dynamics
  simulations using a stochastic reduced-order basis},
\newblock \bibinfo{journal}{Computer Methods in Applied Mechanics and
  Engineering} \bibinfo{volume}{354} (\bibinfo{year}{2019})
  \bibinfo{pages}{37--55}.
\bibitem[{Chernatynskiy et~al.(2013)Chernatynskiy, Phillpot, and
  LeSar}]{chernatynskiy2013uncertainty}
\bibinfo{author}{A.~Chernatynskiy}, \bibinfo{author}{S.~R. Phillpot},
  \bibinfo{author}{R.~LeSar},
\newblock \bibinfo{title}{Uncertainty quantification in multiscale simulation
  of materials: A prospective},
\newblock \bibinfo{journal}{Annual Review of Materials Research}
  \bibinfo{volume}{43} (\bibinfo{year}{2013}) \bibinfo{pages}{157--182}.
\bibitem[{Wang and McDowell(2020)}]{wang2020uncertainty}
\bibinfo{author}{Y.~Wang}, \bibinfo{author}{D.~L. McDowell},
  \bibinfo{title}{Uncertainty Quantification in Multiscale Materials Modeling},
  \bibinfo{publisher}{Woodhead Publishing Limited}, \bibinfo{year}{2020}.
\bibitem[{Falk and Langer(1998)}]{Falk1998}
\bibinfo{author}{M.~L. Falk}, \bibinfo{author}{J.~S. Langer},
\newblock \bibinfo{title}{Dynamics of viscoplastic deformation in amorphous
  solids},
\newblock \bibinfo{journal}{Physical Review E} \bibinfo{volume}{57}
  (\bibinfo{year}{1998}) \bibinfo{pages}{7192--7205}.
\bibitem[{Bouchbinder and Langer(2009)}]{Bouchbinder2009c}
\bibinfo{author}{E.~Bouchbinder}, \bibinfo{author}{J.~S. Langer},
\newblock \bibinfo{title}{Nonequilibrium thermodynamics of driven amorphous
  materials. {III}. shear- transformation-zone plasticity},
\newblock \bibinfo{journal}{Phys. Rev. E} \bibinfo{volume}{80}
  (\bibinfo{year}{2009}) \bibinfo{pages}{031133}.
\bibitem[{Absil et~al.(2004)Absil, Mahony, and Sepulchre}]{Absil2004}
\bibinfo{author}{P.-A. Absil}, \bibinfo{author}{R.~Mahony},
  \bibinfo{author}{R.~Sepulchre},
\newblock \bibinfo{title}{{R}iemannian geometry of {G}rassmann manifolds with a
  view on algorithmic computation},
\newblock \bibinfo{journal}{Acta Applicandae Mathematica} \bibinfo{volume}{80}
  (\bibinfo{year}{2004}) \bibinfo{pages}{199--220}.
\bibitem[{Edelman et~al.(1998)Edelman, Arias, and Smith}]{Edelman1998}
\bibinfo{author}{A.~Edelman}, \bibinfo{author}{T.~A. Arias},
  \bibinfo{author}{S.~T. Smith},
\newblock \bibinfo{title}{The geometry of algorithms with orthogonality
  constraints},
\newblock \bibinfo{journal}{{SIAM}. {J}. {M}atrix {A}nal. \& Appl.}
  \bibinfo{volume}{20} (\bibinfo{year}{1998}) \bibinfo{pages}{303--353}.
\bibitem[{Begelfor and Werman(2006)}]{Begelfor2006}
\bibinfo{author}{E.~Begelfor}, \bibinfo{author}{M.~Werman},
\newblock \bibinfo{title}{Affine invariance revisited},
\newblock \bibinfo{journal}{2006 IEEE Computer Society Conference on Computer
  Vision and Pattern Recognition (CVPR'06)} \bibinfo{volume}{2}
  (\bibinfo{year}{2006}) \bibinfo{pages}{2087--2094}.
\bibitem[{Hamm and Lee(2008)}]{Hamm:2008}
\bibinfo{author}{J.~Hamm}, \bibinfo{author}{D.~D. Lee},
\newblock \bibinfo{title}{Grassmann discriminant analysis: A unifying view on
  subspace-based learning},
\newblock in: \bibinfo{booktitle}{Proceedings of the 25th International
  Conference on Machine Learning}, \bibinfo{publisher}{ACM},
  \bibinfo{year}{2008}, pp. \bibinfo{pages}{376--383}.
\bibitem[{Ye and Lim(2016)}]{YeLim2014}
\bibinfo{author}{K.~Ye}, \bibinfo{author}{L.~H. Lim},
\newblock \bibinfo{title}{{S}chubert varieties and distances between subspaces
  of different dimensions},
\newblock \bibinfo{journal}{{SIAM} {J}ournal on {M}atrix {A}nalysis and
  {A}pplications} \bibinfo{volume}{37} (\bibinfo{year}{2016})
  \bibinfo{pages}{1176--1197}.
\bibitem[{Marrinan et~al.(2014)Marrinan, Beveridge, Draper, and
  Peterson}]{Marrinan2014}
\bibinfo{author}{T.~Marrinan}, \bibinfo{author}{J.~Beveridge},
  \bibinfo{author}{B.~Draper}, \bibinfo{author}{C.~Peterson},
\newblock \bibinfo{title}{Finding the subspace mean or median to fit your
  need},
\newblock \bibinfo{journal}{In CVPR}  (\bibinfo{year}{2014}).
\bibitem[{Karcher(2008)}]{Karcher1977}
\bibinfo{author}{H.~Karcher},
\newblock \bibinfo{title}{Riemannian center of mass and mollifier smoothing},
\newblock \bibinfo{journal}{Communications on Pure and Applied Mathematics}
  \bibinfo{volume}{30} (\bibinfo{year}{2008}) \bibinfo{pages}{509--541}.
\bibitem[{Turaga et~al.(2012)Turaga, Veeraraghavan, Srivastava, and
  Chellappa}]{Turaga2012CorrectionSC}
\bibinfo{author}{P.~K. Turaga}, \bibinfo{author}{A.~Veeraraghavan},
  \bibinfo{author}{A.~Srivastava}, \bibinfo{author}{R.~Chellappa},
\newblock \bibinfo{title}{Correction: Statistical computations on grassmann and
  stiefel manifolds for image and video-based recognition}.
\bibitem[{Von~Luxburg(2007)}]{Luxburg2007}
\bibinfo{author}{U.~Von~Luxburg},
\newblock \bibinfo{title}{A tutorial on spectral clustering},
\newblock \bibinfo{journal}{Statistics and computing} \bibinfo{volume}{17}
  (\bibinfo{year}{2007}) \bibinfo{pages}{395--416}.
\bibitem[{Shi and Malik(2000)}]{Shi2000}
\bibinfo{author}{J.~Shi}, \bibinfo{author}{J.~Malik},
\newblock \bibinfo{title}{Normalized cuts and image segmentation},
\newblock \bibinfo{journal}{IEEE Transactions on Pattern Analysis and Machine
  Intelligence} \bibinfo{volume}{22} (\bibinfo{year}{2000})
  \bibinfo{pages}{888--905}.
\bibitem[{Mohar(1991)}]{Mohar1991}
\bibinfo{author}{B.~Mohar},
\newblock \bibinfo{title}{The laplacian spectrum of graphs},
\newblock \bibinfo{journal}{In Graph theory, combinatorics, and applications}
  \bibinfo{volume}{2} (\bibinfo{year}{1991}) \bibinfo{pages}{871--898}.
\bibitem[{Mohar(1997)}]{Mohar1997}
\bibinfo{author}{B.~Mohar},
\newblock \bibinfo{title}{Some applications of laplace eigenvalues of graphs},
\newblock \bibinfo{journal}{Graph Symmetry: Algebraic Methods and Applications}
  \bibinfo{volume}{497} (\bibinfo{year}{1997}) \bibinfo{pages}{225--275}.
\bibitem[{Ester et~al.(1996)Ester, Kriegel, Sander, and Xu}]{DBSCAN1996}
\bibinfo{author}{M.~Ester}, \bibinfo{author}{H.~P. Kriegel},
  \bibinfo{author}{J.~Sander}, \bibinfo{author}{X.~Xu},
\newblock \bibinfo{title}{A density-based algorithm for discovering clusters in
  large spatial databases with noise},
\newblock \bibinfo{journal}{In: Proceedings of the 2nd International Conference
  on Knowledge Discovery and Data Mining}  (\bibinfo{year}{1996})
  \bibinfo{pages}{226--231}.
\bibitem[{Orszag and Bissonnette(1967)}]{KO1967}
\bibinfo{author}{S.~A. Orszag}, \bibinfo{author}{L.~Bissonnette},
\newblock \bibinfo{title}{Dynamical properties of truncated wiener-hermite
  expansions},
\newblock \bibinfo{journal}{Physics of Fluids} \bibinfo{volume}{10}
  (\bibinfo{year}{1967}) \bibinfo{pages}{2603--2613}.
\bibitem[{Langer(2015)}]{Langer2008}
\bibinfo{author}{J.~S. Langer},
\newblock \bibinfo{title}{Shear-transformation-zone theory of plastic
  deformation near the glass transition},
\newblock \bibinfo{journal}{Physical Review E - Statistical, Nonlinear, and
  Soft Matter Physics} \bibinfo{volume}{77} (\bibinfo{year}{2015})
  \bibinfo{pages}{136--166}.
\bibitem[{Falk and Langer(2011)}]{Falk2011}
\bibinfo{author}{M.~L. Falk}, \bibinfo{author}{J.~S. Langer},
\newblock \bibinfo{title}{Deformation and failure of amorphous, solidlike
  materials},
\newblock \bibinfo{journal}{Annu. Rev. Condens. Matter Phys.}
  \bibinfo{volume}{2} (\bibinfo{year}{2011}) \bibinfo{pages}{353}.
\bibitem[{Bouchbinder and Langer(2009{\natexlab{a}})}]{Bouchbinder2009}
\bibinfo{author}{E.~Bouchbinder}, \bibinfo{author}{J.~S. Langer},
\newblock \bibinfo{title}{Nonequilibrium thermodynamics of driven amorphous
  materials. {I}. internal degrees of freedom and volume deformation},
\newblock \bibinfo{journal}{Phys. Rev. E} \bibinfo{volume}{80}
  (\bibinfo{year}{2009}{\natexlab{a}}) \bibinfo{pages}{031131}.
\bibitem[{Bouchbinder and Langer(2009{\natexlab{b}})}]{Bouchbinder2009b}
\bibinfo{author}{E.~Bouchbinder}, \bibinfo{author}{J.~S. Langer},
\newblock \bibinfo{title}{Nonequilibrium thermodynamics of driven amorphous
  materials. {II}. effective-temperature theory},
\newblock \bibinfo{journal}{Phys. Rev. E} \bibinfo{volume}{80}
  (\bibinfo{year}{2009}{\natexlab{b}}) \bibinfo{pages}{031132}.
\bibitem[{Rycroft et~al.(2008)Rycroft, Sui, and Bouchbinder}]{Rycroft2015}
\bibinfo{author}{C.~H. Rycroft}, \bibinfo{author}{Y.~Sui},
  \bibinfo{author}{E.~Bouchbinder},
\newblock \bibinfo{title}{An eulerian projection method for quasi-static
  elastoplasticity},
\newblock \bibinfo{journal}{Journal of Computational Physics}
  \bibinfo{volume}{30} (\bibinfo{year}{2008}) \bibinfo{pages}{1--14}.
\bibitem[{Boffi and Rycroft(2020)}]{boffi2020parallel}
\bibinfo{author}{N.~M. Boffi}, \bibinfo{author}{C.~H. Rycroft},
\newblock \bibinfo{title}{Parallel three-dimensional simulations of
  quasi-static elastoplastic solids},
\newblock \bibinfo{journal}{Computer Physics Communications}
  (\bibinfo{year}{2020}) \bibinfo{pages}{107254}.
\bibitem[{Hinkle et~al.(2017)Hinkle, Rycroft, Shields, and M.L.}]{Rycroft2017}
\bibinfo{author}{A.~Hinkle}, \bibinfo{author}{C.~Rycroft},
  \bibinfo{author}{M.~Shields}, \bibinfo{author}{F.~M.L.},
\newblock \bibinfo{title}{Coarse graining atomistic simulations of plastically
  deforming amorphous solids},
\newblock \bibinfo{journal}{Phys. Rev. E} \bibinfo{volume}{95}
  (\bibinfo{year}{2017}) \bibinfo{pages}{053001}.
\bibitem[{Shi et~al.(2007)Shi, Katz, Li, and Falk}]{Shi2007}
\bibinfo{author}{Y.~Shi}, \bibinfo{author}{M.~B. Katz},
  \bibinfo{author}{H.~Li}, \bibinfo{author}{M.~L. Falk},
\newblock \bibinfo{title}{{Evaluation of the Disorder Temperature and
  Free-Volume Formalisms via Simulations of Shear Banding in Amorphous
  Solids}},
\newblock \bibinfo{journal}{Physical Review Letters} \bibinfo{volume}{98}
  (\bibinfo{year}{2007}) \bibinfo{pages}{185505}.

\end{thebibliography}







\end{document}